\documentclass[10pt]{amsart}
\usepackage[all]{xy}
\usepackage{graphicx}
\usepackage[colorlinks=true, urlcolor=blue,bookmarks=true,
bookmarksopen=true, citecolor=blue,hypertex]{hyperref}

\usepackage{amssymb,amsmath,amsthm,amscd}

\newcommand{\N}{\mathbb{N}}
\newcommand{\Z}{\mathbb{Z}}
\newcommand{\R}{\mathbb{R}}
\newtheorem{theo}{Theorem}[section]
\newtheorem{cor}[theo]{Corollary}
\newtheorem{prop}[theo]{Proposition}
\newtheorem{rem}[theo]{Remark}

\newtheorem{conj}[theo]{Conjecture}

\begin{document}

\title[Morse theory, Serre spectral sequence, Lagrangian Intersections.]{Homotopical dynamics in symplectic topology.}

\author[J.-F. Barraud,\ O. Cornea]{Jean-Fran\c{c}ois Barraud and Octav Cornea}
\date{May 20, 2005.}

\address{J-F.B.: UFR de Math\'ematiques\newline
\indent Universit\'e de Lille 1\newline \indent 59655 Villeneuve
d'Ascq\newline \indent France} \email{barraud@agat.univ-lille1.fr}

\address{O.C.: University of Montr\'eal\newline
\indent Department of Mathematics and Statistics
\newline \indent CP 6128 Succ. Centre Ville
\newline \indent Montr\'eal, QC H3C 3J7
\newline\indent Canada}
\email{cornea@dms.umontreal.ca}

\begin{abstract}\hskip5pt
This is mainly a survey of recent work on algebraic ways to ``measure''
moduli spaces of connecting trajectories in Morse and Floer theories
as well as related applications to symplectic topology. The paper also
contains some new results. In particular, we show that the methods
of \cite{BaCo} continue to work in general symplectic manifolds (without
any connectivity conditions) but under the bubbling threshold.
\end{abstract}

\maketitle \tableofcontents
\section{Introduction}
The main purpose of this paper is to survey a number of
Morse-theoretic results which show how to estimate algebraically
the high-dimensional moduli spaces of Morse flow lines and to
describe some of their recent applications to symplectic topology.
We also deduce some new applications.

The paper starts with a brief discussion of the various proofs
showing that the differential in the Morse complex is indeed a
differential. With this occasion we introduce the main concepts in
Morse theory and, in particular, the notion of connecting manifold
(or, equivalently, the moduli space of flow lines connecting two
critical points) which is the main object of interest in our
further constructions. Moreover, an extension of one of these
proofs leads naturally to an important result of John Franks
\cite{Fr} which describes the framed cobordism class of connecting
manifolds between consecutive critical points as a certain
relative attaching map. After describing Franks' result, we
proceed to a stronger result initially proved in \cite{Co1} which
computes a framed bordism class naturally associated to the same
connecting manifolds in terms of certain Hopf invariants. While
these results only apply to consecutive critical points we then
describe a recent method to estimate general connecting manifolds
by means of the Serre spectral sequence of the path-loop fibration
having as base the ambient manifold \cite{BaCo}. Some interesting
topological consequences of these results are briefly mentioned as
well as some other methods used in the study of these problems.

The third section discusses a number of symplectic applications.
We start with some results which first appeared in \cite{Co2}.
These use the non-vanishing of certain Hopf invariants
 to deduce the existence of bounded orbits of hamiltonian flows
 (obviously, inside non-compact manifolds).
 This is a very ``soft" type of result even if difficult to prove.
 We then continue in \S\ref{subsec:strips} by describing
 how to use the Serre spectral
 sequence result to detect pseudo-holomorphic strips as well as some
 consequences of the existence of the strips. Most of
 the results of this part have first appeared in \cite{BaCo}
 but there are some that are new: we discuss explicitly the
 detection of pseudoholomorphic strips passing through some
 submanifold and we present a way to construct in a coherent
 fashion our theory for lagrangians in general symplectic
 manifolds as long as we remain {\em under
 a bubbling threshold}. Notice that even the analogue of the classical
 Floer theory (which is a very particular case of our construction)
 has not been explicited  in the  literature
 in the Lagrangian case even if all the necessary ideas
 are present in some form - see \cite{Sch1}
 for the hamiltonian case.

 The paper contains a number of open problems
 and ends with a conjecture which is supported by the
 results in \S\ref{subsec:strips} as well as by recent joint
 results of the second author with Fran\c{c}ois Lalonde.

\section{Elements of Morse theory}\label{sec:Morse}
Assume that $M$ is a compact, smooth manifold without boundary of
dimension $n$. Let $f:M\to \R$ be a smooth Morse function and let
$\gamma: M\times \R\to M$ be a negative gradient Morse-Smale flow
associated of $f$. In particular, $f$ is strictly decreasing along
any non-constant flow line of $\gamma$ and the stable manifolds
$$W^{s}(P)=\{x\in M : \lim_{t\to \infty}\gamma_{t}(x)=P\}$$
and the unstable manifolds
$$W^{u}(Q)=\{x\in M : \lim_{t\to -\infty}\gamma_{t}(x)=Q\}$$
of any pair of critical points $P$ and $Q$ of $f$ intersect
transversally. One of the most useful and simple tools that can be
defined in this context is the Morse complex
$$C(\gamma)=(\Z/2<Crit(f)>,d)~.~$$ Here $\Z/2< S >$ is the
$\Z/2$-vector space generated by the set $S$, the vector space
$\Z/2<Crit(f)>$ has a natural grading given by $|P|=ind_{f}(P),
\forall P\in Crit(f)$ and $d$ is the differential of the complex
which is defined by $$dx=\sum_{|y|=|x|-1}a^{x}_{y}y$$ so that the
coefficients $a^{x}_{y}=\#((W^{u}(x)\cap W^{s}(y))/\R)$. This definition
makes sense because the set $W^{u}(P)\cap W^{s}(Q)$ which consists of all
the points situated on some flow line joining $P$ to $Q$ is invariant by
the $\R$-action given by the flow. Moreover, $W^{u}(P)$ and $W^{s}(P)$
are homeomorphic to open disks which implies that the set
$M^{P}_{Q}=(W^{u}(P)\cap W^{s}(Q))/\R$ has the structure of a smooth (in
general, non-compact) manifold of dimension $|P|-|Q|-1$. We call this
space the moduli space of flow lines joining $P$ to $Q$. It is not
difficult to understand the reasons for the non-compactness of
$M^{P}_{Q}$ when $M$ is compact as in our setting: this is simply due to
the fact that a family of flow lines joining  $P$ to $Q$ may approach a
third, intermediate, critical point $R$. For this to happen it is
necessary (and sufficient - see Smale \cite{Sm} or Franks \cite{Fr}) to
have some flow line which joins $P$ to $R$ and some other joining $R$ to
$Q$. This implies that when $|P|=|Q|+1$ the set $M^{P}_{Q}$ is compact
and thus the sum above is finite.  For further use let's define also the
unstable sphere of a critical point $P$ as $S^{u}_{a}(P)=W^{u}(P)\cap
f^{-1}(a)$ as well as the stable sphere $S^{s}_{a}(P)=W^{s}(P)\cap
f^{-1}(a)$ where $a$ is a regular value of $f$. It should be noted that
this names are slightly abusive as these two sets are spheres, in
general, only if $a$ is sufficiently close to $f(P)$. In that case
$S^{u}_{a}(P)$ is homeomorphic to a sphere of dimension $|P|-1$ and
$S^{s}_{a}(P)$ is homeomorphic to a sphere of dimension $n-|P|-1$. With
this notations ths moduli space $M^{P}_{Q}$ is homeomorphic to
$S^{u}_{a}(P)\cap S^{s}_{a}(P)$ for any $a\in (f(Q),f(P))$ which is a
regular value of $f$.

The main properties of the object defined above are that:
$$d^{2}=0 \ {\rm and}\  H_{\ast}(C(\gamma))\approx
H_{\ast}(M;\Z/2)~.~$$

We will sometimes denote this complex by $C(f)$ and will call it
the (classical) Morse complex of $f$. The flow $\gamma$ may be in
fact even a pseudo- (negative) gradient flow of $f$. There also
exists a version of this complex over $\Z$ in which the counting
of the elements in $M^{x}_{y}$ takes into account orientations.

There are essentially four methods to prove these properties:
\begin{itemize}

\item[(i)] Deducing the equation $\sum_{y}a^{x}_{y}a^{y}_{z}=0$
(which is equivalent to $d^{2}=0$) from the properties of the
moduli spaces $M^{x}_{z}$ with $|x|-|z|=2$.

\item[(ii)] Comparing $a^{x}_{y}$ with a certain relative
attaching map.
\item[(iii)] Expressing $a^{x}_{y}$ in terms of a
connection map in Conley's index theory. \item[(iv)] A method
based on a deformation of the de Rham complex (clearly, in this
case the coefficients are required to be in $\R$).
\end{itemize}

For the rest of this paper the two approaches that are of the most
interest are (i) and (ii). Therefore we shall first say a few
words on the other two methods and will then describe in more
detail the first two. Method (iii) consists in regarding two
critical points $x$, $y$ so that $|x|=|y|+1$ as an
attractor-repellor pair and to apply the general Conley theory of
Morse decompositions to this situation \cite{Sal1}. Method (iv)
has been introduced by Witten in \cite{Wit} and is based on a
deformation of the differential of the de Rham complex which
provides a new differential with respect to which the harmonic
forms are in bijection with the critical points of $f$.

\

Method (i) has been probably folklore for a long time but it first
appeared explicitly in Witten's paper. It is based on noticing
that the moduli space $M^{P}_{Q}$ admits a compactification
$\overline{M}^{P}_{Q}$ which is a compact, topological manifold
with boundary so that the boundary verifies the formula:
\begin{equation}\label{eq:first_bdry}
\partial\overline{M}^{P}_{Q}=\bigcup_{R}\overline{M}^{P}_{R}\times\overline{M}^{R}_{Q}~.~
\end{equation}
There are two main ways to prove this formula. One is analytic and
regards a flow line from $P$ to $Q$ as a solution of a
differential equation $\dot{x}=-\nabla f(x)$ and studies the
properties of such soultions (this method is presented in the book
of Schwarz \cite{Sch}). A second approach is more
topological/dynamical in nature as is described in detail by Weber
\cite{We}. Clearly, from formula (\ref{eq:first_bdry}) we
immediately deduce $\sum_{y}a^{x}_{y}a^{y}_{z}=0$ and hence
$d^{2}=0$. Just a little more work is needed to deduce from here
the second property.

\

Method (ii) was the one best known classically and it is
essentially implicit in Milnor's $h$-cobordism book \cite{Mil}. It
is based on the observation that $a^{x}_{y}$ can be viewed as
follows. First, to simplify slightly the argument assume that the
only critical points in  $f^{-1}([f(y),f(x)])$ are $x$ and $y$. It
is well known that for $a\in (f(y),f(x))$, there exists a
deformation retract
$$r:M(a)=f^{-1}(-\infty,a]\to M(f(y)-\epsilon)\cup_{\phi_{y}} D^{|y|}=M'$$ where the attaching map
\begin{equation}\label{eq:attaching}
\phi_{y}: S^{u}_{f(y)-\epsilon}(y)\to M(f(y)-\epsilon)
\end{equation} is just the inclusion and
$\epsilon$ is small. This deformation retract follows the flow
till reaching of $U(W^{u}(y))\cup M(f(y)-\epsilon)$ where
$U(W^{u}(y)$ is a tubular neighbourhood of $W^{u}(y)$ so that the
flow is transverse to its boundary and then collapses this
neighbourhood to $W^{u}(y)$ by the canonical projection. Clearly,
applying this remark to each critical point of $f$ provides a
$CW$-complex of the same homotopy type as that of $M$ and with one
cell $\bar{x}$ for each element of $x\in Crit(f)$. To this
cellular decomposition we may associate a celullar complex
$(C'(f), d')$ with the property that $d\bar{x}=\sum
k^{\bar{x}}_{\bar{y}}\bar{y}$ where $k^{\bar{x}}_{\bar{y}}$ is, by
definition, the degree of the composition:
\begin{equation}\label{eq:rel_attaching}
\psi^{x}_{y}:S^{u}_{a}(x)\stackrel{\phi_{x}}{\longrightarrow}
M(a)\stackrel{r}{\longrightarrow}M' \stackrel{u}{\longrightarrow}
M'/M(f(y)-\epsilon) \end{equation} with $\phi_{x}:S^{u}_{a}(x)\to
M(a)$ again the inclusion and where the last map, $u$, is the
projection onto the respective topological quotient space (which
is homeomorphic to the sphere $S^{|y|}$).

Notice now that $M^{x}_{y}\subset S^{u}_{a}(x)$ is a finite union
of points say $P_{1},\ldots, P_{k}$. Imagine a small disk
$D_{i}\subset S^{u}_{a}(x)$ around $P_{i}$. The key (but
geometrically clear) remark is that the composition of the flow
$\gamma$ together  with the retraction $r$ transports $D_{i}$ (if
it is chosen sufficiently small) homeomorphically  onto a
neighbourhood  of $y$ inside $W^{u}(y)$. Therefore, the degree of
$deg(\psi^{x}_{y})=a^{x}_{y}$ and thus $d=d'$ which shows that $d$
is a differential and that the homology it computes agrees with
the homology of $M$.

\

As we shall see further, the points of view reflected in the
approaches at (i) and (ii) lead to interesting applications which
go much beyond ``classical" Morse theory. Method (iv), while
striking and inspiring appears for now not to have been exploited
efficiently.

\subsection{Connecting Manifolds}
One way to look to the Morse complex is by viewing the
coefficients $a^{x}_{y}$ of the differential as a measure of the
$0$-dimensional manifold $M^{x}_{y}$. The question we discuss here
is in what way we can measure algebraically the similar higher
dimensional moduli spaces. This is clearly a significant issue
because, obviously, only a very superficial part of the dynamics
of the negative gradient flow of $f$ is encoded in the
$0$-dimensional moduli spaces of connecting flow lines.

\

As a matter of terminology, the space $M^{P}_{Q}$ when viewed
inside the unstable sphere $S^{u}_{a}(P)$ (with $f(P)-a$ positive
and very small) is also called \cite{Fr} the connecting manifold
of $P$ and $Q$.

It was mentioned above that, in general, a connecting manifold
$M^{P}_{Q}$ is not closed. However, if the critical points $P$ and
$Q$ are {\em consecutive} in the sense that there does not exist a
critical point $R$ so that $M^{P}_{R}\times
M^{R}_{Q}\not=\emptyset$, then $M^{P}_{Q}$ is closed.

\subsubsection{Framed Cobordism Classes}

An important remark of John Franks \cite{Fr} is that connecting
manifolds are canonically framed. First recall that a framed
manifold $V$ is a submanifold $V\hookrightarrow S^{n}$ which has a
trivial normal bundle together with a trivialization of this
bundle. Two such trivializations are equivalent (and will
generally be identified) if they are restrictions of a
trivialization of the normal bundle of $V\times [0,1]$ inside
$S^{n}\times [0,1]$. We also recall the Thom-Pontryagin
construction in this context \cite{Mil2}. Assuming
$V\hookrightarrow S^{n}$ is framed we define a map
$$\phi_{V}:S^{n}\to S^{codim(V)}$$
as follows: consider a tubular neighbourhood $U(V)$ of $V$, use
the framing to define a homeomorphism $\psi:U(V)\to D^{k}\times V$
where $D^{k}$ is the closed disk of dimension $k=codim(V)$,
consider the composition
$\psi':U(V)\stackrel{\psi}{\longrightarrow}D^{k}\times
V\stackrel{p_{1}}{\longrightarrow} D^{k}\to D^{k}/S^{k-1}=S^{k}$
and define $\phi_{V}$ by extending $\psi'$ outside $U(V)$ by
sending each $x\in S^{n}\backslash U(V)$ to the base point in
$D^{k}/S^{k-1}=S^{k}$. The homotopy class of this map is the same
if two framings are equivalent. It is easy to see that two framed
manifolds (of the same dimension) are cobordant iff their
associated Thom maps are homotopic.

\

We return now to Franks' remark and notice that the manifolds
$M^{P}_{Q}$ are framed. First, we make the convention to view
$M^{P}_{Q}$ as a submanifold of the unstable sphere of $P$,
$S^{u}_{a}(P)$ (the other choice would have been to use $S^{s}(Q)$
as ambient manifold). Notice that we have
$$M^{P}_{Q}=S^{u}_{a}(P)\cap W^{s}(Q)$$ and this intersection is
transversal. Clearly, as $W^{s}(Q)$ is homeomorphic to a disk, its
normal bundle in $M$ is trivial and any two trivializations of
this bundle are equivalent. This implies that the normal bundle of
$M^{P}_{Q}\hookrightarrow S^{u}_{a}(P)$ is also trivial and a
trivialization of the normal bundle of $W^{s}(Q)$ provides a
trivialization of this bundle which is unique up to equivalence.

\

Recall that if $P$ and $Q$ are consecutive critical points, then
$M^{P}_{Q}$ is closed. As we have seen that it is also framed we
may associate to it a framed cobordism class
$$\widetilde{M^{P}_{Q}}\in \pi_{|P|-1}(S^{|Q|})~.~$$
Moreover, it is easy to see that the function $f$ may be perturbed
without modifying the dynamics of the the negative gradient flow
so that the cell attachments corresponding to the critical points
$Q$ and $P$ are in succession. Therefore the map
$$\psi^{P}_{Q}:S^{|P|-1}\to S^{|Q|}$$ defined as in formula
(\ref{eq:rel_attaching}) is still defined. The main result of
Franks in \cite{Fr} is:

\begin{theo}\label{theo:fr}\cite{Fr} Assume $P$ and $Q$ are consecutive critical
points of $f$. Up to sign $\widetilde{M^{P}_{Q}}$ coincides with
the homotopy class of $\psi^{P}_{Q}$.
\end{theo}

The idea of proof of this result is quite simple. All that is required is
to make even more precise the constructions used in the approach (ii)
used to show $d^{2}=0$ for the Morse complex. For this we fix for
$W^{s}(Q)$ a normal framing $o$ which is invariant by translation along
the flow $\gamma$ and which at $Q\in W^{s}(Q)$ is given by a basis $e$ of
$T_{Q}W^{u}(Q)$ (this is possible because $W^{u}(Q)$ and $W^{s}(Q)$
intersect transversally at $Q$). We also fix the tubular neighbourhood
$U(W^{u}(Q))$ so that the projection $r'':U(W^{u}(Q))\to W^{u}(Q)$ has
the property that $(r'')^{-1}(Q)=W^{s}(Q)\cap U(W^{u}(Q))$ and, for any
point $y\in (r'')^{-1}(Q)$, we have $(r'')_{\ast}(o_{y})=e$. Moreover, we
may assume that the normal bundle of $M^{P}_{Q}$ in $S^{u}_{a}(P)$ is
just the restriction of the normal bundle of $W^{s}(Q)$ (in fact, the two
are, in general, only isomorphic and not equal but this is just a minor
issue).

Now, follow what happens with the framing of $M^{P}_{Q}$ along the
composition $u\circ r$. For this we write $r=r''\circ r'$ where $r'$
follows the flow till reaching $U(W^{u}(Q))$. Now pick a point in $x\in
M^{P}_{Q}$ together with its normal frame $o_{x}$ at $x$. After applying
$r'$, the pair $(x,o_{x})$ is taken to a pair $(x',o_{x'})$ with
$x'\in\partial ((r'')^{-1}(Q))$. Applying now $r''$, the image of
$(x',o_{x'})$ is $(Q,e)$. Take now $V$ a tubular neighbourhood of
$M^{P}_{Q}\hookrightarrow S^{u}_{a}(P)$ together with an identification
$V\approx D^{|Q|}\times M^{P}_{Q}$ which is provided by the framing $o$.
The argument above implies that if the constant $\epsilon$ used to
construct the map $u: M(f(Q)-\epsilon)\cup W^{u}(Q)=M'\to
M'/M(f(Q)-\epsilon)$ is very small, then the composition $u\circ r''\circ
r'$ coincides with the relevant Thom-Pontryagin map.

\subsubsection{Framed Bordism Classes and Hopf invariants}

It it natural to wonder whether, besides their framing, there are
some other properties of the connecting manifolds which can be
detected algebraically.

A useful point of view in this respect turns out to be the
following: imagine the elements of $M^{P}_{Q}$ as path or loops on
$M$. The fact that they are paths is obvious (we parametrize them
by the value of $-f$; the negative sign gives the flow lines the
orientation coherent with the negative gradient) but they can be
transformed into loops very easily. Indeed, fix a simple path in
$M$ which joins all the critical points of $f$ and contract this
to a point thus obtaining a quotient space $\widehat{M}$ which has
the same homotopy type as $M$.  Let $q:M\to \widehat{M}$ be the
quotient map. We denote by $\Omega M$ the space of based loops on
$M$ and keep the notation $q$ for the induced map $\Omega M\to
\Omega\widehat{M}$. This discussion shows that there are
continuous maps
$$j^{P}_{Q}:M^{P}_{Q}\to \Omega\widehat{M}~.~$$
These maps have first been defined and used in \cite{Co1} and they
have some interesting properties.

For example, given such a map $j^{P}_{Q}$ and assuming that $P$
and $Q$ are consecutive it is natural to ask whether the homology
class $[M^{P}_{Q}]\in H_{|P|-|Q|-1}(\Omega M; \Z)$ is
computable (here $[M^{P}_{Q}]$ is the image by $j^{P}_{Q}$ of the
 fundamental class of $M^{P}_{Q}$). We shall see that quite a bit more is indeed
possible: the full framed bordism class associated to $j^{P}_{Q}$
and to the canonical framing on $M^{P}_{Q}$ can be expressed as a
relative Hopf invariant.

\

To explain this result first recall that if $V\hookrightarrow
S^{n}$ is framed and $l:V\to X$ is a continuous map with $V$ a
closed manifold we may construct a richer Thom-Pontryagin map as
follows. We again consider a tubular neighbourhood $U(V)$ of $V$
in $S^{n}$ together with an identification $U(V)\approx
D^{k}\times V$ where $k=codim(V)$ which is provided by the
framing. We now define a map
$$
\bar{\phi}_{V}:U(V)\stackrel{j}{\longrightarrow} D^{k}\times V
\stackrel{id\times l}{\longrightarrow} D^{k}\times X\to
(D^{k}\times X)/(S^{k-1}\times X)$$ where $j:U(V)=D^{k}\times
V\to V$ is the projection and the last map is just the quotient map
(which identifies $S^{k-1}\times X$ to the base point). Notice that
$\bar{\phi}_{V}(\partial U(V))\subset S^{k-1}\times X$. Therefore, we may
extend the definition above to a map
$$\bar{\phi}_{V}:S^{n}\to (D^{k}\times X)/(S^{k-1}\times X)$$ by sending
all the points in the complement of $U(V)$ to the base point. It is
well-known (and a simple exercise of elementary homotopy theory) that
there exists a (canonical) homotopy equivalence $(D^{k}\times
X)/(S^{k-1}\times X)\simeq \Sigma^{k}(X^{+})$ where $\Sigma X$ is the
(reduced) suspension of $X$, $\Sigma^{i}$ is the suspension iterated
$i$-times and $X^{+}$ is the space $X$ with an added disjoint point
(notice also that $\Sigma^{k}(X^{+})=\Sigma^{k}X\vee S^{k}$ where $\vee$
denotes the {\em wedge} or the one point union of spaces). This allows us
to view the map $\bar{\phi}_{V}$ as a map with values in
$\Sigma^{k}(X^{+})$. The framed bordism class of $V$ is simply the
homotopy class $[\bar{\phi}_{V}]\in \pi_{n}(\Sigma^{k}(X^{+}))$. This is
independent of the various choices made in the construction. Two pairs of
data (framings included) $(V,l)$ and $(V',l')$ admit an extension to a
manifold $W\subset S^{n}\times [0,1]$ with $\partial W=V\times\{0\}\cup
V'\times\{1\}$ iff $\bar{\phi}_{V}\simeq \bar{\phi}_{V'}$. Notice also
that to an element $\alpha\in \pi_{n}(\Sigma^{k} X^{+})$ we may associate
a homology class $[\alpha]\in H_{n-k}(X)$ obtained by applying
the Hurewicz homomorphism, desuspending $k$-times and projecting
on the $H_{\ast}(X)$ term in $H_{\ast}(X^{+})$.

\

Returning now to our connecting manifolds $M^{P}_{Q}$ we again
focus on the case when $P$ and $Q$ are consecutive. The map
$j^{P}_{Q}$ together with the canonical framing provide a homotopy
class $\{M^{P}_{Q}\}\in \pi_{|P|-1}(\Sigma^{|Q|}(\Omega M ^{+}))$
(to simplify notation we have replaced $\widehat{M}$ with $M$ here
- the two are homotopy equivalent).

\

As indicated above, it turns out that this class can be computed
in terms of a relative Hopf invariant. We shall now discuss how
this invariant is defined.

Assume that $S^{q-1}\stackrel{\alpha}{\to} X_{0}\to X'$ and
$S^{p-1}\stackrel{\beta}{\to} X'\to X''$ are two successive cell
attachments and that $X''$ is a subspace of some larger space $X$.
In particular,  $X'= X_{0}\cup_{\alpha} D^{q}$, $X''= X'\cup_{\beta}
D^{p}$. Let $S\subset D^{q}$ be the $q-1$-sphere of radius $1/2$.
There is an important map called the coaction associated to
$\alpha$ which is defined by  the composition  $$\nabla: X'\to
X'/S \approx S^{q}\vee X'$$ where the first map identifies all the
points of $S$ to a single one and the second is a homeomorphism
(in practice it is convenient to also assume that all the maps and
spaces involved are pointed and in that case we view $D^{q}$ as
the reduced cone over $S^{q-1}$).

We consider the composition
$$\psi(\beta,\alpha):S^{p-1}\stackrel{\beta}{\longrightarrow}
X'\stackrel{\nabla}{\longrightarrow}S^{q}\vee X'\stackrel{id\vee
j}{\longrightarrow} S^{q}\vee X''\hookrightarrow S^{q}\vee X$$ and
notice that if $p_{2}:S^{q}\vee X\to X$ is the projection on the
second factor, then the composition $p_{2}\circ
\psi(\beta,\alpha)$ is null-homotopic. This is due to the fact
that this composition is homotopic to
$S^{p-1}\stackrel{\beta}{\longrightarrow}X'\to X''\hookrightarrow
X$ which is clearly null-homotopic. We now consider the map
$p_{2}:S^{q}\vee X\to X$. It is well-known in homotopy theory that
any map may be transformed into a fibration. In our case this
comes down to considering the free path fibration
$$t: \widetilde{P}X\to X$$ where $\widetilde{P}X$ is the set of
all continuous path in $X$ parametrized by $[0,1]$, $t(\gamma)=\gamma(0)$.
We take the pull back of this
fibration over $p_{2}$. The total space $\widetilde{E}$ of the
resulting fibration has the same homotopy type as $S^{q}\vee X$
and it is endowed with a canonical projection
$$\widetilde{p}:\widetilde{E}\to X$$ which replaces $p_{2}$, $\widetilde{p}(z,\gamma)=\gamma(1)$.
It is an  exercise in homotopy theory to see that the fibre of the
fibration $\widetilde{p}$ is homotopic to $\Sigma^{q}((\Omega
X)^{+})$ and that, moreover, the inclusion of this
fibre in the total space is injective in homotopy. As the composition
$p_{2}\circ\psi(\beta,\alpha)$ is homotopically trivial, the
homotopy exact sequence of the fibration $\widetilde{E}\to X$ implies
that $\psi(\beta,\alpha)$ admits a lift to
$\bar{\psi}(\alpha,\beta):S^{p-1}\to \Sigma^{q}((\Omega X)^{+})$
whose homotopy class does not depend on the choice of lift. We let
$H(\alpha,\beta)\in \pi_{p-1}(\Sigma^{q}((\Omega X)^{+}))$  be
equal to this homotopy class and we call it the relative Hopf
invariant associated to the attaching maps $\alpha$ and $\beta$
(for a discussion of the relations between this Hopf invariants
and other variants see  Chapters 6 and 7 in \cite{CLOT}).

\

To return to Morse theory, recall from (\ref{eq:attaching}) that
passing through the two consecutive critical points $Q$ and $P$
leads to two successive attaching maps $\phi_{Q}:S^{|Q|-1}\to
M(f(Q)-\epsilon)$ and $\phi_{P}:S^{|P|-1}\to M(f(P)-\epsilon)$ (we
assume again - as we may - that the set $f^{-1}([f(Q),f(P)])$ does
not contain any other critical points besides $P$ and $Q$).
Moreover, as we know the inclusion $M'=M(f(Q)-\epsilon)\cup
W^{u}(Q)\hookrightarrow M(f(P)-\epsilon)$ is a homotopy
equivalence. Therefore, the construction above can be applied to
$\phi_{Q}$ and $\phi_{P}$ and it leads to a relative Hopf
invariant $$H(P,Q)\in \pi_{|P|-1}(\Sigma ^{|Q|}(\Omega M^{+}))~.~$$

With these constructions our statement is:

\begin{theo}\label{theo:bordism}\cite{Co1}\cite{Co2}
The homotopy class $H(P,Q)$ coincides (up to sign) with the
bordism class $\{M^{P}_{Q}\}$. In particular, the homology class
$[M^{P}_{Q}]$ equals (up to sign and desuspension) the Hurewicz
image of $H(P,Q)$.
\end{theo}

The proof of this result can be found in \cite{Co2} (a variant
proved by a slightly different method appears in \cite{Co1}). The
proof is considerably more complicated than the proof of Theorem
\ref{theo:fr} so we will only present a rough justification here.
To simplify notation we let $M_{0}=M(f(Q)-\epsilon)$. Let
$M_{1}=M(f(Q)-\epsilon)\cup U(W^{u}(Q))$.
 Recall, that the inclusions $M'\hookrightarrow M_{1}
 \hookrightarrow M(f(P)-\epsilon)$ are homotopy equivalences.

Let $\mathcal{P}: \Omega M\to PM\to M$ be the path-loop fibration
(of total space the {\em based} paths on $M$ and of fibre the
space of based loops on $M$). We denote by $E_{0}$ the total space
of the pull-back of the fibration $\mathcal{P}$ over the inclusion
$M_{0}\subset M$. Similarly, we let $E_{1}$ be the total space of
the pull-back of $\mathcal{P}$ over the inclusion $M_{1}\to M$.
The key remark is that the attaching map $\phi_{P}: S^{u}(P)\to
M_{1}$ admits a natural lift to a map
$\widetilde{\phi_{P}}:S^{u}(P)\to E_{1}$. Indeed, we assume here
that all the critical points are identified to the base point. The
space $E_{1}$ consists of the based paths in $M$ that end at
points in $M_{1}$. But each element of the image of $\phi_{P}$
corresponds to precisely such a path which is explicitly given by
the corresponding flow line (we need to use here Moore paths and
loops which are paths parametrized by arbitrary intervals $[0,a]$
and not only the interval $[0,1]$). Consider the inclusion
$E_{0}\hookrightarrow E_{1}$. It is not difficult to see that the
quotient topological space $E_{1}/E_{0}$ admits a canonical
homotopy equivalence $\eta: E_{1}/E_{0}\to \Sigma^{|Q|}(\Omega M
^{+})$. Therefore, we may consider the composition
$\eta'=\eta\circ\widetilde{\phi_{P}}$. It is possible to show that
this map $\eta'$ is homotopic to $H(P,Q)$. At the same time, we
see that the restriction of $\widetilde{\phi_{P}}$ to $M^{P}_{Q}$
coincides with $j^{P}_{Q}$. Moreover, by making explicit $\eta$ it
is also possible to see that $\widetilde{\phi_{P}}$ is homotopic
to the Thom-Pontryagin map associated to $M^{P}_{Q}$.

\subsubsection{Some topological applications}
\label{subsubsec:top_appl} We now decribe a couple  of topological
applications of Theorem \ref{theo:bordism}. The idea behind both
of them is quite simple: the function $-f$ is also a Morse
function and the critical points $Q$, $P$ are conscutive critical
points for $-f$. Therefore, the connecting manifold $M^{Q}_{P}$ is
well defined as well as its associated bordism homotopy class
$\{M^{Q}_{P}\}$. Clearly, the underlying space for both
$M^{P}_{Q}$ and $M^{Q}_{P}$ is the same. The map $j^{Q}_{P}$ is
different from $j^{P}_{Q}$ just by reversing the direction of the
loops. The two relevant framings may also be different. The
relation between them is somewhat less straightforward but it
still may be understood by considering $M$ embedded inside a high
dimensional euclidean space and taking into account the twisting
induced by the stable normal bundle. In all cases, this
establishes a relationship between the two Hopf invariants
$H(P,Q)$ and $H(Q,P)$.

\

 A. The first application \cite{Co1} concerns the construction
of examples of non-smoothable, simply-connected, Poincar\'e
duality spaces. The idea is as follows: we construct Poincar\'e
duality spaces which have a simple $CW$-decomposition and with the
property that for certain two successively attached cells $e,f$
the resulting Hopf invariant $H$ and the Hopf invariant $H'$
associated to the dual cells $f',e'$ are not related in the way
described above. If the respective Poincar\'e duality space is
smoothable, then the given cell decomposition can be viewed as
associated to an appropriate Morse function and this leads to a
contradiction. The obstructions to smoothability constructed in
this way are obstructions to the lifting of the Spivak normal
bundle to $BO$. This is an obstruction theory problem but one
which can be very difficult to solve directly in the presence of
many cells. Thus, this approach is quite efficient to construct
examples.

\

B. The second application \cite{Co3} concerns the detection of
obstructions to the embedding of $CW$-complexes in euclidean
spaces in low codimension. The argument in this case goes roughly
as follows. If the $CW$-complex $X$ embedds in $S^{n}$, then we
may consider a neighbourhood $U(X)$ of $X$ which is a smooth
manifold with boundary. We consider a smooth Morse function
$f:U(X)\to \R$ which is constant, maximal and regular on the
boundary of $U(X)$. If $P$ and $Q$ are two consecutive critical
points for this funtion we obtain that
$\Sigma^{k}H(P,Q)=\Sigma^{k'}H(Q,P)$ for certain values of $k$ and
$k'$ which can be estimated explicitly - the main reason for this
equality is that the Morse function in question is defined on the
sphere so all the questions involving the stable normal bundle
become irrelevant. If $X$ admits some reasonably explicit
cell-decomposition it is possible to express $H(P,Q)$ as the Hopf
invariant $H$ of some successive attachment of two cells $e,f$ and
$H(Q,P)$ as $\Sigma^{k''}H'$ where $H'$ is another similar Hopf
invariant. The obstructions to embedding appear because the low
codimension condition translates to the fact that $k'+k''> k$.
This can be viewed as an obstruction because it means that after
$k$ suspensions the homotopy class of $H$ has to de-suspend more
than $k$-times.

\begin{rem} {\rm The applications at A and B are purely of homotopical type.
It is natural to expect that the Morse theoretical arguments that
were used to establish these statements can be replaced by purely
homotopical ones but this has not been done till now.}
\end{rem}

\subsubsection{The Serre spectral sequence}\label{subsubsec:serre_ss}

Theorem \ref{theo:bordism} provides considerable information on
connecting manifolds for pairs of {\em consecutive} critical
points. However, it does not shed any light on the case of
non-consecutive ones. Clearly, if the critical points are not
consecutive the respective connecting manifold is not closed and
thus no bordism or cobordism class can be directly associated to
it. However, after compactification, the boundary of this
connecting manifold has a special structure reflected by equation
(\ref{eq:first_bdry}). As we shall see following \cite{BaCo}, this
structure is sufficient to construct an algebraic invariant which
provides an efficient ``measure" of all connecting manifolds.

\

This construction is based on the fact that the maps
$$j^{P}_{Q}:M^{P}_{Q}\to \Omega M$$ are
compatible with compactification and with the formula
(\ref{eq:first_bdry}) in the following sense. Recall that here
$\Omega M$ are the based Moore loops on $M$ (these are loops
parametrized by intervals $[0,a]$), the critical points of $f$
have been identified to a single point and, moreover, in the
definition of $j^{P}_{Q}$ we use the parametrization of the flow
lines by the values of $-f$. Recall that we have a product given by
the concatenation of loops
$$\mu:\Omega M\times \Omega M\to \Omega M~.~$$

With these notations it is easy to see that we have the following
formula:
\begin{equation}\label{eq:product}
j^{P}_{Q}(u,v)=\mu(j^{P}_{R}(u),j^{R}_{Q}(v))
\end{equation}
where $(u,v)\in \overline{M}^{P}_{R}\times
\overline{M}^{R}_{Q}\subset\partial\overline{M}^{P}_{Q}$.

\

We proceed with our consruction. Let $C_{\ast}(X)$ be the
(reduced) cubical complex of $X$ with coefficients in $\Z/2$.
Notice that there is a natural map $$C_{k}(X)\otimes C_{k'}(Y)\to
C_{k+k'}(X\times Y)~.~$$

\

A family of cubical chain $s^{x}_{y}\in
C_{|x|-|y|-1|}(\overline{M}^{x}_{y})$, $x,y\in Crit(f)$ is called
a {\em representing chain system} for the moduli spaces
$M^{x}_{y}$ if for each pair of critical points $x,z$ we have:
\begin{itemize}
\item[i.]
$$d s^{x}_{z}=\sum_{y} s^{x}_{y}\otimes s^{y}_{z}$$
\item[ii.] $s^{x}_{z}$ represents the fundamental class in
$H_{|x|-|z|-1}(M^{x}_{z},\partial M^{x}_{z})$. \end{itemize}

\

It is easy to show by induction on the index difference $|x|-|z|$
that such representing chain systems exist. We now fix such a
representing chain system $\{s^{x}_{y}\}$ and we define
$a^{x}_{y}\in C_{|x|-|y|-1}(\Omega M)$ by
$$a^{x}_{y}=(j^{x}_{y})_{\ast}(s^{x}_{y})~.~$$

Notice that this definition extends the definition of these
coefficients in the usual Morse case when $|x|-|y|-1=0$. We have a
product map
$$\cdot:C_{k}(\Omega M)\otimes C_{k'}(\Omega M)\longrightarrow
C_{k+k'}(\Omega M\times \Omega
M)\stackrel{C_{\ast}(\mu)}{\longrightarrow}C_{k+k'}(\Omega M)$$
which makes $C_{\ast}(\Omega M)$ into a differential ring. The
discussion above shows that inside this ring we have the formula
$$da^{x}_{z}=\sum_{y}a^{x}_{z}\cdot a^{z}_{y}~.~$$
An elegant way to rephrase this formula is to group these
coefficients in a matrix $A=(a^{x}_{y})$ and then we have
\begin{equation}\label{eq:diff_matrix}
dA=A^{2}~.~\end{equation}

We now define a new chain complex $\mathcal{C}(f)$ associated to
$f$ by
\begin{equation}\label{eq:ext_cplx}\mathcal{C}(f)=(C_{\ast}(\Omega
M)\otimes \Z/2<Crit(f)>\ ,d)\ \ , \ \ dx=\sum_{y}a^{x}_{y}\otimes
y~.~
\end{equation}

We shall call this complex {\em the extended Morse complex of $f$}. Here,
$C_{\ast}(\Omega M)\otimes \Z/2<Crit(f)>$ is viewed as a graded
$C_{\ast}(\Omega M)$-module and $d$ respects this structure in the sense
that it verifies $d(a\otimes x)=(da)\otimes x + a (dx)$ (the grading on
$Crit(f)$ is given, as before, by the Morse index). Choosing orientations
on all the stable manifolds of all the critical points induces a
co-orientation on all the unstable manifolds, and hence an orientations
on the intersections $W^{u}(P)\cap W^{s}(Q)$ and finally on all the
moduli spaces $M^{P}_{Q}$~: we may then use $\Z$-coefficients for this
complex as well as, of course, for the classical Morse complex. In this
case appropriate signs appear in the formulae above. Clearly, $d^{2}=0$
due to (\ref{eq:diff_matrix}).

\

By definition, the coefficients $a^{x}_{y}$ represent the moduli
spaces $M^{x}_{y}$. However, these coefficients are not invariant
with respect to the choices made in their construction. Therefore,
it is remarkable that there is a natural construction which
extracts from this complex a useful algebraic invariant which is
not just the homology of the complex - as it happens, this homology is
not too interesting as it coincides with that of a point.

Consider the obvious differential filtration which is defined on
this complex by $$F^{k}\mathcal{C}(f)=C_{\ast}(\Omega M)\otimes
\Z/2< x\in Crit(f) :\ ind_{f}(x)\leq k>~.~$$ Denote the associated
spectral sequence by
$$E(f)=(E^{r}_{p,q}(f),d^{r})~.~$$

\begin{theo}\label{theo:serre_ss}\cite{BaCo} When $M$ is simply connected and if
 $r\geq 2$ the spectral sequence $E(f)$ coincides with the
 Serre spectral sequence of the path-loop fibration
 $$\mathcal{P}:\Omega M\to PM\to M~.~$$
\end{theo}

\begin{rem} \label{rem:serre_ss} {\rm a. A similar result can be
established even in the absence of the simple-connectivity
condition which has been assumed here to avoid some technical
complications.

b. The Serre spectral sequence of the path-loop fibration of a
space $X$ contains considerable information on the homotopy type
of the space. In particular, there are spaces with the same
cohomology and cup-product but which may be distinguished by their
respective Serre spectral sequences. }\end{rem}

To outline the proof of the theorem we start by recalling the
construction of the Serre spectral sequence in the form which will
be of use here. We shall assume here that the Morse function $f$
is self-indexed (in the sense that for each critical point $x$ we
have $ind_{f}(x)=k \Rightarrow f(x)=k$) and that it has a single
minimum denoted by $m$. Let $M_{k}=f^{-1}((-\infty, k+\epsilon])$.
We have $$M_{k}=M_{k-1}\bigcup_{\phi_{y}} D^{k}_{y}$$ where the
union is taken over all the critical points $y\in Crit_{k}(f)$ and
$\phi_{y}:S^{u}(y)\to M_{k-1}$ are the respective attaching maps. Denote
by $E_{k}$ the total space of the fibration induced by pull-back over the
inclusion $M_{k}\hookrightarrow M$ from the fibration $\mathcal{P}$.
Consider the filtration of $C_{\ast}(PM)$ given by
$F^{k}P=Im(C_{\ast}(E_{k})\longrightarrow C_{\ast}(PM))$. The spectral
sequence induced by this filtration is
 invariant after the second page and is
 precisely the Serre spectral sequence (this spectral sequence may be
 constructed as above but by using an arbitrary skeletal filtration $\{X_{k}\}$ of
 a space $X$ which has the same homotopy type as that of $M$; in our case
 the particular filtration given by the sets $M_{k}$ is a natural choice).

\

For further use, we also notice that there is an obvious action of
$\Omega M$ on $PM$ and this action induces one on each $E_{k}$.
Therefore, we may view $C_{\ast}(E_{k})$ as a $C_{\ast}(\Omega
M)$-module.

\

The first step in proving the theorem is to consider a certain
compactification of the unstable manifolds of the critical points
of $f$. Recall that $f$ is self-indexed and that $m$ is the unique
minimum critical point of $f$. Fix $x\in Crit(f)$ and define the
following equivalence relation on the set
$\overline{M}^{x}_{m}\times [0,f(x)]$:

$$(a,t)\sim (a',t') \ {\rm iff}\ t=t' \ {\rm and}\
a(-\tau)=a'(-\tau)\ \forall \tau\geq t~.~$$ Here the elements of
$\overline{M}^{x}_{m}$ are viewed as paths in $M$ parametrized by
the value of $-f$ (so that $f(a(-\tau))=\tau$).

Denote by $\widehat{W}(x)$ the resulting quotient topological
space. Notice that, if $y\in \overline{W^{u}(x)}$, then there
exists some $a\in \overline{M}^{x}_{m}$ so that $y$ is on the
(possibly broken) flow line represented by $a$. Or, in other
words, so that $a(-f(y))=y$. This path $a$ might not be unique.
Indeed, inside $\widehat{W}(x)$ there is precisely one equivalence
class $[a,f(y)]$ (with $a(-f(y))=y$) for each (possibly) broken
flow line joining $x$ to $y$. Clearly, if $y\in W^{u}(x)$, then
there is just one such flow line and so the natural surjection
$$\pi: \widehat{W}(x)\longrightarrow \overline{W^{u}(x)},\
\pi([a,t])=a(-t)$$ is a homeomorphism when restricted to
 $\pi^{-1}(W^{u}(x))$. Thus we may view $\widehat{W}(x)$ as a
special compactification of $W^{u}(x)$ or as a desingularization
of $\overline{W^{u}(x)}$.  It is not difficult to believe (but
harder to show and will not be proven here see \cite{BaCo}) that
$\widehat{W}(x)$ is a topological manifold with boundary and
moreover
\begin{equation}\label{eq:bdry_cpct}
\partial \widehat{W}(x)=\bigcup_{y}M^{x}_{y}\times \widehat{W}(y)~.~
\end{equation}

\

We continue the proof of the Theorem \ref{theo:serre_ss} with the
remark that there are obvious maps
$$h_{x}: \widehat{W}(x)\to PM$$
which associate to $[a,t]$ the path in $M$ which follows $a$ from
$x$ to $a(-t)$. These maps and the maps $j^{x}_{y}$ are compatible
with formula (\ref{eq:bdry_cpct}) in the sense that
\begin{equation}\label{eq:compat}
h_{x}(a',[a'',t])=j^{x}_{y}(a')\cdot h_{y}([a'',t])
\end{equation}
where $(a',[a'',t])\in M^{x}_{y}\times \widehat{W}(y)$ and $\cdot$
represents the action of $\Omega M$ on $PM$.

Now, we may obviously rewrite (\ref{eq:bdry_cpct}) as:
$$\partial \widehat{W}(x)=\bigcup_{y} \overline{M}^{x}_{y}\times \widehat{W}(y)~.~$$

Given the representing chain system $\{s^{x}_{y}\}$ it is easy to
construct an associated representing chain system for
$\widehat{W}(x)$. This is a system of chains $v(x)\in
C_{|x|}(\widehat{W}(x))$ so that $v(x)$ represents the fundamental
class of $C_{|x|}(\widehat{W}(x),\partial\widehat{W}(x))$ and we
have the formula
$$dv(x)=\sum_{y}s^{x}_{y}\otimes v(y)~.~$$

 Finally, we define
a $C_{\ast}(\Omega M)$-module chain map
$$\alpha : \mathcal{C}(f)\to C_{\ast}(PM)$$
 by $$\alpha(x)=(h_{x})_{\ast}(v(x))~.~$$

 The formulas above show that we have
 $$d[(h_{x})_{\ast}(v(x))]=\sum_{y}a^{x}_{y}\cdot (h_{y})_{\ast}(v(y))$$
and so $\alpha$ is a chain map. It is clear that the map $\alpha$
is filtration preserving and it is not difficult to see that it
induces an isomorphism at the $E^{2}$ level of the induced
spectral sequences and this concludes the proof of Theorem
\ref{theo:serre_ss}.

\begin{rem} {\rm
a. Another important but immediate property of $\widehat{W}(x)$ is
that it is a contractible space. Indeed, all the points in
$\overline{M}^{x}_{m}\times \{f(x)\}$ are in the same equivalence
class. Moreover, each point $[a,t]\in\widehat{W}(x)$ has the
property that it is related by the path $[a,\tau], \ \tau\in
[t,f(x)]$ to $\ast=[a,f(x)]$. The contraction of $\widehat{W}(x)$
to $\ast$ is obtained by deforming $\widehat{W}(x)$ along these
paths. Given that $\widehat{W}(x)$ is a contractible topological
manifold with boundary, it is natural to suspect that
$\widehat{W}(x)$ is homeomorphic to a disk. This is indeed the
case as is shown in \cite{BaCo} and is an interesting fact in
itself because it implies that the union of the unstable manifolds
of a self-indexed Morse-Smale function gives a $CW$-decomposition
of $M$. The attaching map of the cell $\widehat{W}(x)$ is simply
the restriction of $\pi$ to $\partial \widehat{W}(x)$.

b. The Serre spectral sequence result above and the bordism result
in Theorem \ref{theo:bordism} are obviously related via the
central role of the maps $j^{P}_{Q}$. There is also a more
explicit relation. Indeed, (a stable version of) the Hopf
invariants appearing in Theorem \ref{theo:bordism} can be
interpreted as differentials in the Atyiah-Hirzebruch-Serre
spectral sequence of the path-loop fibration with coefficients in
the stable homotopy of $\Omega M$. Moreover, the relation
(\ref{eq:first_bdry}) can be understood as also keeping track of
the framings. This leads to a type of extended Morse complex in
which the coefficients of the differential are stable Hopf
invariants \cite{Co1}. All of this strongly suggests that the
construction of the complex $\mathcal{C}(f)$ can be enriched so as
to include the framings of the connecting manifolds and, by the
same method as above, that the whole Atyiah-Hirzebruch-Serre
spectral sequence should be recovered from this construction.

c. Another interesting question, open even for consecutive
critical points $P$, $Q$, is whether there are some additional
constraints on the topology of the connecting manifolds
$M^{P}_{Q}$ besides those imposed by Theorem \ref{theo:bordism}.

d. Yet another open question is how this machinery can be adapted
to the Morse-Bott situation and how it can be extended to general
Morse-Smale flows (not only gradient-like ones).

e. It is natural to wonder what is the richest level of
information that one can extract out of the moduli spaces of Morse
flow lines. At a naive level, the union of all the points situated
on the flow lines of $f$ is precisely the whole underlying
manifold $M$ so we expect that there should exist some assembly
process producing the manifold $M$ out of these moduli spaces.
Such a machine has been constructed by Cohen, Jones and Segal
\cite{CoJoSe1}\cite{CoJoSe2}. They show that one can form a
category out of the moduli spaces of connecting trajectories and
that the classifying space of this category is of the
homeomorphism type of the underlying manifold. In their
construction an essential point is that the glueing of flow lines
is associative. This approach is quite different from the
techniques above and does not imply the results concerning the
extended Morse complex or the Hopf invariants that we have
presented. The two points of view are, essentially, complementary.
}\end{rem}

To end this section it is useful to make explicit a relation
between Theorems \ref{theo:serre_ss} and \ref{theo:bordism} (we
assume as above that $M$ is simply-connected).

\begin{prop}\label{prop:hlgy} Assume that there are $q,p\in \N$ so that
 in the Serre spectral sequence of the path loop fibration of
 $M$ we have $E^{2}_{k,s}=0$ for $q<k<p$ and there is an element
 $a\in E^{2}_{p,0}$ so that $d^{p-q}a\not=0$, then any Morse-Smale
 function on $M$ has a pair of consecutive critical points $P$, $Q$
 of indexes at least $q$ and at most $p$ so that the homology class
 $[M^{P}_{Q}]\in H_{|P|-|Q|-1}(\Omega M)\not=0$.
\end{prop}

Clearly, Theorem \ref{theo:serre_ss} directly implies that, even
without any restriction on $E^{2}$, if we have $d^{r}a\not=0$ with
$a\in E^{r}_{p,0}$, then for any Morse-Smale function $f$ there
are critical points $P$ and $Q$ with $|P|=p$ and $|Q|=p-r$ so that
$M^{P}_{Q}\not=\emptyset$. Indeed, if this would not be the case, then all
the coefficients $a^{x}_{y}$ in the extended Morse complex of $f$
are null whenever $|x|=p$, $|y|=p-r$. By the construction of the
associated spectral sequence, this leads to a contradiction.
However, the pair $P$, $Q$ resulting from this argument might have
a connecting manifold which is not closed so that its homology
class is not even defined and, thus, Proposition \ref{prop:hlgy}
provides a stronger conclusion. The proof of the Proposition is as
follows. Recall that $E^{2}_{s,r}\approx H_{s}(M)\otimes
H_{r}(\Omega M)$ and so $H_{\ast}(M)=0$ for $q\leq 0\leq p$. If
there are some points $P, Q\in Crit(f)$ with $q\leq |Q|,|P|\leq p$
so that the differential of $P$ in the classical Morse complex
contains $Q$ with a non-trivial coefficient then this pair $P$,
$Q$ may be taken as the one we are looking for. If all such
differentials in the classical Morse complex are trivial it
follows that the critical points of index $p$ and $q$ are
consecutive. In this case, the geometric arguments used in the proofs of
either Theorem \ref{theo:bordism} or \ref{theo:serre_ss} imply
that if for all pairs $P$, $Q$ with $|P|=p$, $|Q|=q$ we would have
$[M^{P}_{Q}]=0$, then the differential $d^{p-q}$ would vanish on
$E^{p-q}_{p,0}$.

\begin{rem} \label{rem:consec} {\rm Notice that the pair of critical points $P$ and $Q$
constructed in the proposition verify the property that $|P|$ and
$|Q|$ are consecutive inside the set $\{ ind_{f}(x): x\in
Crit(f)\}$.}
\end{rem}

\subsection{Operations} We discuss here a different and, probably, more familiar
approach to understanding connecting manifolds as well as other
related Morse theoretic moduli spaces. This point of view has been
used extensively by many authors - Fukaya \cite{Fu}, Betz and
Cohen \cite{BeCo} being just a few of them. For this reason we
shall review this technique very briefly.

\

Given two consecutive critical points $x$, $y$ notice that the set
$T^{x}_{y}=W^{u}(x)\cap W^{s}(y)$ is homeomorphic to the
un-reduced suspension of $M^{x}_{y}$. Therefore, we may see this
as an obvious inclusion
$$i^{x}_{y}:\Sigma M^{x}_{y}\to M$$ and we may consider the homology class
$[T^{x}_{y}]=(i^{x}_{y})_{\ast}(s[M^{x}_{y}])$ where $s$ is
suspension and $[M^{x}_{y}]$ is the fundamental class. There
exists an obvious evaluation map
$$e:\Sigma \Omega M\to M$$ which is induced by
$\Omega M\times [0,1]\to M,\ (\beta,t)\to \beta(t)$ (the loops
here are parametrized by the interval $[0,1]$ but this is a minor
technical difficulty). It is easy to see, by the definition of
this evaluation map, that
$[T^{x}_{y}]=e_{\ast}((j^{x}_{y})_{\ast}([M^{x}_{y}]))$. In
general the map $e_{\ast}$ is not injective in homology. Clearly,
the full bordism class $\{M^{x}_{y}\}$ carries much more
information than the homology class $[T^{x}_{y}]$. Still, there is
a direct way to determine $[T^{x}_{y}]$ without passing through a
calculation of $\{M^{x}_{y}\}$ and we will now describe it.

Consider a second Morse-Smale function $g:M\to \R$ so that its
associated unstable and stable manifolds $W^{u}_{g}(-)$,
$W^{s}_{g}(-)$ intersect transversally the stable and unstable
manifolds of $f$ and, except if they are of top dimension, they
avoid the critical points of $f$.

\

Fix $x,y\in Crit(f)$ and $s\in Crit(g)$ so that
$|x|-|y|-ind_{g}(s)=0$. We may define $k(x,y;s)=\#(T^{x}_{y}\cap
W^{s}_{g}(s)$  (where the counting takes into account the relevant
orientations if we work over $\Z$). We now put
$$\bar{k}^{x}_{y}=\sum_{s} k(x,y;s)s\in C(g)~.~$$

The essentially obvious claim is that:
\begin{prop}\label{prop:hlgy_class} The chain $\bar{k}^{x}_{y}$
is a cycle whose homology class is $[T^{x}_{y}]$.
\end{prop}

Indeed we have $\sum_{s}k(x,y;s)h^{s}_{z}=0$ where
$ind_{g}(z)=ind_{g}(s)-1$ and $h^{s}_{z}$ are the coefficients in
the classical Morse complex of $g$. This equality is valid because
we may consider the $1$-dimensional space $T^{x}_{y}\cap
W^{s}_{g}(z)$. This is an open $1$-dimensional manifold whose
compactification is a $1$-manifold whose boundary points are
counted precisely by the formula $\sum_{s}k(x,y;s)h^{s}_{z}=0$. By
basic intersection theory it is immediate to see that the homology
class represented by this cycle is $[T^{x}_{y}]$.

\

While this construction does not shed a lot of light on the
properties of $M^{x}_{y}$ its role is important once we use it to
recover the various homological operations of $M$. To see how this
is done from our perspective notice that the intersection
$$T^{x}_{y}\cap W^{s}_{g}(s)=W^{u}_{f}(x)\cap W^{s}_{f}(y)\cap
W^{s}_{g}(s)$$ can be viewed as a particular case of the following
situation: assume that $f_{1}$, $f_{2}$, $f_{3}$ are three
Morse-Smale functions in general position and define
$$T^{x,y}_{z}=(W^{u}_{f_{1}}(x)\cap W^{u}_{f_{2}}(y)\cap
W^{s}_{f_{3}}(z))~.~$$  If we assume that $|x|+|y|-|z|-n=0$ we may
again count the  points in $T^{x,y}_{z}$ with appropriate signs
and we may define coefficients $t^{x,y}_{z}=\#T^{x,y}_{z}$.

This leads to an operation \cite{BeCo}\cite{Fu}
$$C(f_{1})\otimes C(f_{2})\to C(f_{3})$$
given as a linear extension of $$x\otimes y\to \sum_{z}
t^{x,y}_{z}z~.~$$ It is easy to see that this operation descends
in homology and that it is in fact the dual of the $\cup$-product.
Moreover, the space $T^{x,y}_{z}$ may be viewed as obtained by
considering a graph formed by three oriented edges meeting into a
point with the first two entering the point and the other one
exiting it and considering all the configurations obtained by
mapping this graph into $M$ so that to the first edge we associate
a flow line of $f_{1}$ which exits $x$, to the second edge a flow
line of $f_{2}$ which exits $y$ and to the third a flow line of
$f_{3}$ which enters $z$. Clearly, this idea may be pushed further
by considering other, more complicated graphs and understanding
what are the operations they correspond to as was done by Betz and
Cohen \cite{BeCo}.

\section{Applications to Symplectic topology}
We start with some applications that are rather ``soft" even if
difficult to prove and we shall continue in the main part of the
section,\S\ref{subsec:strips}, with some others that go deeper.
\subsection{Bounded orbits}
 We fix a symplectic
manifold $(M,\omega)$ which is not compact. Assume
that $H:M\to \R$ is a smooth hamiltonian whose associated
hamiltonian vector field is denoted by $X_{H}$. One of the main
questions in hamiltonian dynamics is whether a given regular
hypersurface $A=H^{-1}(a)$ of $H$ has any closed caracteristics,
or equivalently, whether the hamiltonian flow of $H$ has any
periodic orbits in $A$. As $M$ is not compact, from the point of
view of dynamical systems, the first natural question is whether
$X_{H}$ has any {\em bounded} orbits in $A$.  Moreover, there is a
remarkable result of Pugh and Robinson \cite{PuRo}, the
$C^{1}$-closing lemma, which shows that, for a generic choice of
$H$, the presence of bounded orbits insures the existence of some
periodic orbits.  Therefore, we shall focus in this subsection on
the detection of bounded orbits. It should be noted however that
the detection of periodic orbits in this way is not very effective
because the periods of the orbits found can not be estimated.
Moreover, there is no reasonable test to decide whether a given
hamiltonian belongs to the generic family to which the
$C^{1}$-closing lemma applies. Finally, it will be clear from the
methods of proof described below that these results are also soft
in the sense that they are not truly specific to Hamiltonian flows
but rather they apply to many other flows.

\

An example of a bounded orbit result is the following statement
\cite{Co2}.

\begin{theo}\label{theo:bounded}
Assume that $H$ is a Morse-Smale function with respect to a
riemannian metric $g$ on $M$ so that $M$ is metrically complete
and there exists an $\epsilon$ and a compact set $K\subset M$ so
that $||\nabla_{g} H(x)||\geq \epsilon$ for $x\not\in K$. Suppose
that $P$ and $Q$ are two critical points of $H$ so that $|P|$ and
$|Q|$ are successive in the set $\{ind_{H}(x) : x\in Crit(H)\}$.
If the stabilization $[H(P,Q)]\in \pi^{S}_{|P|-|Q|-1}(\Omega M)$
of the Hopf invariant $H(P,Q)$ is not trivial then there are
regular values $v\in (H(Q),H(P))$ so that $H^{-1}(v)$ contains
bounded orbits of $X_{H}$.
\end{theo}

Before describing the proof of this result let's notice that the
theorem is not difficult to apply. Indeed, one simple way to
verify that there are pairs $P$, $Q$ as required is to use
Proposition \ref{prop:hlgy} together with Remark \ref{rem:consec}
with a minor adaptation required in a non-compact setting. This
adaptation consists of replacing the Serre spectral sequence of
the path loop fibration with the Serre spectral sequence of a
relative fibration $\Omega M\to (E_{1},E_{0})\to (N_{1},N_{0})$
where $N_{1}$ is an isolating neighbourhood for the gradient flow
of $H$ which contains $K$ and $N_{0}$ is a (regular) exit set for
this neighbourhood (to see the precise definition of these Conley
index theoretic notions see \cite{Sal1}). The fibration is induced
by pull-back from the path-loop fibration $\Omega M\to PM\to M$
over the inclusion $(N_{1},N_{0})\hookrightarrow (M,M)$. In short,
because the gradient of $H$ is away from $0$ outside of a compact set,
pairs $(N_{1},N_{0})$ as above are easy to produce and if the pair
$(N_{1},N_{0})$ has some interesting topology it is easy to deduce
the existence of non-constant bounded orbits. Here is a concrete
example.

\begin{cor}\label{cor:bounded_orb}
Assume that $M$ is the contangent bundle of some closed,
simply-connected
 manifold $N$ of dimension $k\geq 2$ and $\omega$ is  an arbitrary symplectic form.
 Assume that $H:M\to \R$ is Morse and that outside of some compact set containing
 the $0$ section, $H$ restricts to each fibre of the bundle to a non-degenerate
 quadratic form. Then, $X_{H}$ has bounded, non constant orbits.
\end{cor}

Of course, this result is only interesting when there are  no
compact level hypersurfaces of $H$. This does happen if the
quadratic form in question has an index which is neither $0$ nor
$k$. The proof of this result comes down to the fact that as $N$
is closed and not a point there exists a lowest dimensional
homology class $u\in H_{t}(N)$ which is transgressive in the Serre
spectral sequence ( this means $d^{t}u\not=0$). Using the
structure of function quadratic at infinity of $H$ it is easy to
construct a pair $(N_{1},N_{0})$ where $N_{1}$ is is homotopic to
a disk bundle of base $N$ and $N_{0}$ is the associated sphere
bundle. The spectral sequence associated to this pair can be
related by the Thom isomorphism to the Serre spectral sequence of
the path-loop fibration over $N$ and the element $\bar{u}\in
H_{\ast}(N_{1},N_{0})$ which corresponds to $u$ by the Thom
isomorphism will have a non-vanishing differential. This means
that Proposition \ref{prop:hlgy} may be used to show the non-triviality of
a homology class $[M^{P}_{Q}]$ for $P$ and $Q$ as in
Theorem \ref{theo:bounded}.
By Theorem \ref{theo:bordism}, $[M^{P}_{Q}]$ is the same up to
sign as the homology class of the Hopf invariant $H(P,Q)$ so
Theorem \ref{theo:bounded} is applicable to detect bounded orbits.

\

We now describe the proof of the theorem. The basic idea of the
proof is simpler to present in the particular case when
$H^{-1}(H(Q),H(P))$ does not contain any critical value. In this
case let $A=H^{-1}(a)$ where $a\in (H(Q),H(P))$. We intend to show
that $A$ contains some bounded orbits of $X_{H}$. To do this
notice that the two sets $S_{1}=W^{u}(P)\cap A$,
$S_{2}=W^{s}(Q)\cap A$ are both diffeomorphic to spheres. We now
assume that no bounded orbits exist and we consider a compact
neighbourhood $U$ of $S_{1}\cup S_{2}$. Assume that we let $S_{2}$
move along the flow $X_{H}$. As this flow has no bounded orbits,
each point of $S_{2}$ will leave $U$ at some moment. Suppose that
we are able to perturb the flow induced by $X_{H}$ to a new
deformation $\eta : M\times \R\to M$ so that for some finite time
$T$ {\em all} the points in $S_{2}$ are taken simultaneously
outside $U$ (in other words $\eta_{T}(S_{2})\cap U=\emptyset$) and
so that $\eta$ leaves $Q$ fixed. It is easy to see that this
implies that $S_{1}\cap S_{2}$ is bordant to the empty set which,
by Theorem \ref{theo:bordism}, is impossible because
$H(P,Q)\not=0$. This perturbation $\eta$ is in fact not hard to
construct by using some elements of Conley's index theory and the
fact that the maximal invariant set of $X_{T}$ inside $U$ is the
empty set (the main step here is to possibly modify also $U$ so
that it admits a regular exit region $U_{0}\subset U$ and we then
construct $\eta$ so that it follows the flow lines of $X_{T}$ but
stops when reaching $U_{0}$, this eliminates the problem of
``bouncing" points which first exit $U$ but later re-enter it).

The case when there are critical points in $H^{-1}(H(Q),H(P))$
follows the same idea but is considerably more difficult. The main
difference comes from the fact that the sets $S_{1}$ and $S_{2}$
might not be closed manifolds. Even their closures $\bar{S_{1}}$
and $\bar{S_{2}}$ are not closed manifolds in general but might be
singular sets. To be able to proceed in this case we first replace
$P$ and $Q$ with a pair of critical points of the same index so
that for any critical point $Q'\in H^{-1}(H(Q),H(P))$ with
$ind_{H}(Q')=ind_{H}(Q)$ we have $[H(P,Q')]=0$. We then take $a$
very close to $H(Q)$ so that $S_{2}$ at least is diffeomorphic to
a sphere. We then first study the stratification of $\bar{S_{1}}$:
there is a top stratum of dimension $|P|-1$ which is $S_{1}$ and a
singular stratum $S'$ of dimension $|Q|-1$ which is the union of
the sets $W^{u}(Q')\cap A$ for all $Q'$ so that
$M^{P}_{Q'}\not=\emptyset$ and $|Q'|=|Q|$. Notice that the way to
construct the null-bordism of $S_{1}\cup S_{2}$ is to consider in
$A\times [0,T]$ the submanifold $W=(\eta_{t}(\bar{S_{2}}),t)$ and
intersect it with $W'=\bar{S_{1}}\times [0,T]$  - we assume here
$\eta_{T}(S_{2})\cap S_{1}=\emptyset$. Clearly, we need this
intersection to be transverse and this can be easily achieved by a
perturbation of $\eta$. The main technical difficulty is that $L$
might intersect the singular part, $S'\times [0,T]$. Indeed,
$dim(W)=n-q$, $dim(S')=q-1$, $dim(A)=n-1$ and so generically the
intersection $I$ between $W$ and $S'\times [0,T]$ is
$0$-dimensional and not necessarily void. By studying the geometry
around each of the points of $I$ it can be seen that $S_{1}\cap
S_{2}$ is bordant to the union of the $M^{P}_{Q'}$'s where $Q'\in
H^{-1}(H(Q),H(P))$ (roughly, this follows by eliminating from the
singular bordism $W\cap W'$ a small closed, cone-like
neighbourhood around each singular point and showing that the
boundary of this cone-like neighbourhood is homeomorphic to a
$M^{P}_{Q'}$). We now use the fact that all the stable bordism
classes of the $M^{P}_{Q'}$'s vanish (because $[H(P,Q')]=0$) and
this leads us to a contradiction. Notice also that, at this point,
we need to use stable Hopf invariants (or bordism classes) $\in
\pi^{S}(\Omega M)$ because, by contrast to the stable case, the
unstable Thom-Pontryagin map associated to a disjoint union is not
necessarily equal to the sum of the Thom-Pontryagin maps of the
terms in the union and hence, unstably, even if we know
$H(P,Q')=0, \forall Q'$ we still can not deduce $H(P,Q)=0$.

\begin{rem}{\rm
It would be interesting to see whether, under some additional
assumptions, a condition of the type $[H(P,Q)]\not=0$ implies the
existence of periodic orbits and not only bounded ones.}
\end{rem}

\subsection{Detection of pseudoholomorphic
strips and Hofer's norm}\label{subsec:strips} In this subsection
we shall again use the Morse theoretic techniques described in \S
\ref{sec:Morse} and, in particular, Theorem \ref{theo:serre_ss} to
study some symplectic phenomena by showing that Floer's complex
can be enriched in a way similar to the passage from the classical
Morse complex to the extended one.

\subsubsection{Elements of Floer's theory.}
We start by recalling very briefly some elements from Floer's
construction (for a more complete exposition see, for example,
\cite{Sal2}).

We shall assume from now on that $(M,\omega)$ is a symplectic manifold -
possibly non-compact but in that case convex at infinity - of dimension
$m=2n$. We also assume that $L,L'$ are closed (no boundary, compact)
Lagrangian submanifolds of $M$ which intersect transversally.

To start the description of our applications it is simplest to assume for
now that $L,L'$ are simply-connected and that
$\omega|_{\pi_{2}(M)}=c_{1}|_{\pi_{2}(M)}=0$. Cotangent bundles of
simply-connected manifolds offer immediate examples of manifolds
verifying these conditions.

We fix a path $\eta\in \mathcal{P}(L,L')=\{\gamma\in C^{\infty}([0,1], M)
: \gamma(0)\in L$, $\gamma(1)\in L'\}$ and let $\mathcal{P}_{\eta}(L,L')$
be the path-component of $\mathcal{P}(L,L')$ containing $\eta$. This path
will be trivial homotopically in most cases, in particular, if $L$ is
hamiltonian isotopic to $L'$. We also fix an almost complex structure $J$
on $M$ compatible with $\omega$ in the sense that the bilinear form
$X,Y\to \omega(X,JY)=\alpha (X,Y)$ is a Riemannian metric. The set of all
the almost complex structures on $M$ compatible with $\omega$ will be denoted
by $\mathcal{J}_{\omega}$. Moreover, we also consider a smooth
$1$-periodic Hamiltonian $H:[0,1]\times M\to \R$ which is constant
outside a compact set and its associated $1$-periodic family of
hamiltonian vector fields $X_{H}$ determined by the equation
$$\omega(X^{t}_{H},Y)=-dH_{t}(Y) \ , \ \forall Y~.~$$ we denote by
$\phi_{t}^{H}$ the associated Hamiltonian isotopy. We also assume
that $\phi^{H}_{1}(L)$ intersects transversally $L'$.

In our setting, the action functional below is well-defined:
\begin{equation}\label{eq:action}
\mathcal{A}_{L,L',H}:\mathcal{P}_{\eta}(L,L')\to \R \ , \ x\to
-\int \overline{x}^{\ast}\omega +\int_{0}^{1}H(t,x(t))dt
\end{equation}
where $\overline{x}(s,t):[0,1]\times [0,1]\to M$ is such that
$\overline{x}(0,t)=\eta(t)$, $\overline{x}(1,t)=x(t)$, $\forall
t\in [0,1]$, $x([0,1],0)\subset L$, $x([0,1],1)\subset L'$. The
critical points of $\mathcal{A}$ are the orbits of $X_{H}$ that
start on $L$, end on $L'$ and which belong to
$\mathcal{P}_{\eta}(L,L')$. These orbits are in bijection with a
subset of $\phi_{1}^{H}(L)\cap L'$ so they are finite in number.
If $H$ is constant these orbits coincide with the intersection
points of $L$ and $L'$ which are in the class of $\eta$. We denote
the set of these orbits by $I(L,L'; \eta, H)$ or shorted $I(L,L')$
if $\eta$ and $H$ are not in doubt.

We now consider the solutions $u$ of Floer's equation:
\begin{equation}\label{eq:floer}
\frac{\partial u}{\partial s}+J(u)\frac{\partial u}{\partial
t}+\nabla H(t,u)=0
\end{equation}
with
 $$ u(s,t):\R\times [0,1]\to M \ , u(\R,0)
 \subset L \ , \ u(\R,1)\subset L' ~.~$$
 When
$H$ is constant, these solutions are called
\emph{pseudo-holomorphic strips}.

For any strip  $u\in\mathcal{S}(L,L')=\{u\in C^{\infty}(\R\times
[0,1],M) : u(\R,0)\subset L \ , \ u(\R,1)\subset L' \}$  consider
the energy
\begin{equation}\label{eq:energy}
E_{L,L',H}(u)=\frac{1}{2}\int_{\R\times [0,1]} ||\frac{\partial
u}{\partial s}||^{2}+ ||\frac{\partial u}{\partial
t}-X^{t}_{H}(u)||^{2}\ ds\ dt ~.~
\end{equation}
For a generic choice of $J$, the solutions $u$ of (\ref{eq:floer})
which {\em are of finite energy}, $E_{L,L',H}(u)<\infty$, behave
like negative gradient flow lines of $\mathcal{A}$. In particular,
$\mathcal{A}$ decreases along such solutions. We consider the
moduli space
\begin{equation}\label{eq:param_moduli}
\mathcal{M}'=\{u\in\mathcal{S}(L,L') : u \ {\rm verifies
(\ref{eq:floer}) } \ , \ E_{L,L',H}(u)<\infty \}~.~
\end{equation}
The translation $u(s,t)\to u(s+k,t)$ obviously induces an $\R$
action on $\mathcal{M}'$ and we let $\mathcal{M}$ be the quotient
space. For each $u\in \mathcal{M}'$ there exist $x,y\in
I(L,L';\eta, H)$ such that the (uniform) limits verify
\begin{equation}\label{eq:ends}
\lim_{s\to-\infty}u(s,t)=x(t) \ , \
\lim_{s\to+\infty}u(s,t)=y(t)~.~
\end{equation}
We let $\mathcal{M}'(x,y)=\{u\in \mathcal{M}' : u \ {\rm verifies}
\ (\ref{eq:ends})\}$ and $\mathcal{M}(x,y)=\mathcal{M}'(x,y)/\R$
so that $\mathcal{M}=\bigcup_{x,y} \mathcal{M}(x,y)$. If needed,
to indicate to which pair of Lagrangians, to what Hamiltonian and
to what almost complex structure are associated these moduli
spaces we shall add $L$ and $L'$, $H$, $J$ as subscripts (for
example, we may write $\mathcal{M}_{L,L',H,J}(x,y)$).

\

For $x,y\in I(L,L';\eta,H)$ we let
$$\begin{matrix}\mathcal{S}(x,y)=\{ u\in C^{\infty}([0,1]\times
[0,1], M) \ :& u([0,1],0)\subset L, u([0,1],1)\subset L', \\
\hspace{55pt} u(0,t)=x(t) ,u(1,t)=y(t)\}~.~ &
\end{matrix}$$
To each $u\in S(x,y)$ we may associate its Maslov index
$\mu(u)\in\Z$ \cite{Viterbo} and it can be seen that, in our
setting, this number only depends on the points $x,y$. Thus, we
let $\mu(x,y)=\mu(u)$. Moreover, we have the formula
\begin{equation}\label{eq:transit}
\mu(x,z)=\mu(x,y)+\mu(y,z)~.~
\end{equation}
According to these relations, the choice of an arbitrary intersection
point $x_{0}$ and the normalization $|x_{0}|=0$, defines a grading
$|.|$ such that~:
$$
 \mu(x,y)=|x|-|y|.
$$

There is a notion of regularity for the pairs of $(H,J)$ so that, when
regularity is assumed, the spaces $\mathcal{M}'(x,y)$ are smooth
manifolds (generally non-compact) of dimension $\mu(x,y)$ and in this
case $\mathcal{M}(x,y)$ is also a smooth manifold of dimension
$\mu(x,y)-1$. Regular pairs $(H,J)$ are generic and, in fact, they are so
even if $L$ and $L'$ are not transversal (but in that case $H$ can not be
assumed to be constant), for example, when $L=L'$.

\

Floer's construction is natural in the following sense. Let
$L''=(\phi_{1}^{H})^{-1}(L')$. Consider the map $b_{H}:
\mathcal{P}(L,L'')\to \mathcal{P}(L,L')$ defined by
$(b_{H}(x))(t)=\phi_{t}^{H}(x(t))$. Let
$\eta'\in\mathcal{P}(L,L'')$ be such that $\eta=b_{H}(\eta')$.
Clearly, $b_{H}$ restricts to a map between
$\mathcal{P}_{\eta'}(L,L'')$ and $\mathcal{P}_{\eta}(L,L')$ and it
restricts to a bijection $I(L,L'';\eta',0)\to I(L,L';\eta, H)$.

It is easy to also check
$$\mathcal{A}_{L,L',H}(b_{H}(x))=\mathcal{A}_{L,L'',0}(x)$$ and
that the map $b_{H}$ identifies the geometry of the two action
functionals. Indeed for $u:\R\times [0,1]\to M$ with
$u(\R,0)\subset L$, $u(\R,1)\subset L''$,
$\tilde{u}(s,t)=\phi_{t}(u(s,t))$, $\tilde{J}=\phi^{\ast}J$ we
have
$$\phi_{\ast}(\frac{\partial u}{\partial s}+ \tilde{J} \frac{\partial u}{\partial t}) =
\frac{\partial \tilde{u}}{\partial s}+J (\frac{\partial
\tilde{u}}{\partial t}-X_{H})~.~$$ Therefore, the map $b_{H}$
induces diffeomorphisms:
$$b_{H}:\mathcal{M}_{L,L'',\tilde{J},0}(x,y)\to \mathcal{M}_{L,L',J,H}(x,y)$$
where we have identified $x,y\in L\cap L''$ with their orbits
$\phi^{H}_{t}(x)$ and $\phi^{H}_{t}(y)$.

\

Finally, the non-compactness of $\mathcal{M}(x,y)$ for $x,y\in
I(L,L';\eta,H)$ is due to the fact that, as in the Morse-Smale
case, a sequence of strips $u_{n}\in \mathcal{M}(x,y)$ might
``converge" (in the sense of Gromov) to a broken strip. There are
natural compactifications of the moduli spaces $\mathcal{M}(x,y)$
called Gromov compactifications and denoted by
$\overline{\mathcal{M}}(x,y)$ so that each of the spaces
$\overline{\mathcal{M}}(x,y)$ is a topological manifold with
corners
 whose boundary verifies:
\begin{equation}\label{eq:Grom_comp}
\partial\overline{\mathcal{M}}(x,y)=\bigcup_{z\in I(L,L';\eta,H)}\ \overline{\mathcal{M}}(x,z)
\times \overline{\mathcal{M}}(z,y)~.~
\end{equation}
 A complete proof of this fact can be found in
\cite{BaCo} (when $dim(\mathcal{M}^{x}_{y})=1$ the proof is due to
Floer and is now classical).

\subsubsection{Pseudoholomorphic strips and the Serre spectral sequence}
We will now construct a complex $\mathcal{C}(L,L';H,J)$ by a
method that mirrors the construction of $\mathcal{C}(f)$ in
\S\ref{subsubsec:serre_ss}.

This complex, called the {\em extended Floer complex} associated
to $L,L',H,J$ has the form:
$$\mathcal{C}(L,L';H,J)=(C_{\ast}(\Omega L)\otimes \Z/2<I(L,L';\eta, H)>, D)$$
where the cubical chains $C_{\ast}(\Omega L)$ have, as before,
$\Z/2$-coefficients.
If needed, the moduli spaces $\mathcal{M}(x,y)$ can be endowed with
orientations which are compatible with formula (\ref{eq:Grom_comp}), and
so we could as well use $\Z$-coefficients.

%

\

To define the differential we first fix a simple path $w$ in $L$
which joins all the points $\gamma(0)$, $\gamma\in I(L,L';H)$ and
we identify all these points to a single one by collapsing this
path to a single point. We shall continue to denote the resulting
space by $L$ to simplify notation.

For each moduli space $\mathcal{M}(x,y)$ there is a continuous map
$$l^{x}_{y}:\mathcal{M}(x,y)\to \Omega L$$
which is defined by associating to $u\in \mathcal{M}(x,y)$ the
path $u(\R,0)$ parametrized by the (negative) values of the action
functional $\mathcal{A}$. This is a continuous map and it is seen
to be compatible with formula (\ref{eq:Grom_comp}) in the same
sense as in (\ref{eq:compat}).

We pick a representing chain system $\{k^{x}_{y}\}$ for the moduli
spaces $\mathcal{M}(x,y)$ and we let
$$m^{x}_{y}=(l^{x}_{y})_{\ast}(k^{x}_{y})\in C_{\ast}(\Omega L)$$
and
\begin{equation}\label{eq:differential_coeff}
Dx=\sum_{y}m^{x}_{y}\otimes y~.~
\end{equation}
As in the case of the extended Morse complex the fact that
$D^{2}=0$ is an immediate consequence of formula
(\ref{eq:Grom_comp}).

\begin{rem}\label{rem:coeff}{\rm  a. There is an apparent assymetry
between the roles of $L$ and $L'$ in the definition of this extended Floer complex. In fact,
the coefficients of this complex belong naturally to an even
bigger and more symmetric ring than $C_{\ast}(\Omega L)$. Indeed,
consider the space $T(L,L')$ which is the homotopy-pullback of the
two inclusions $L\hookrightarrow M$, $L'\hookrightarrow M$. This
space is homotopy equivalent to the space of all the continuous
paths $\gamma : [0,1]\to M$ so that $\gamma(0)\in L$,
$\gamma(1)\in L'$. By replacing both $L$ and $M$ by the respective
spaces obtained by contracting the path $w$ to a point, we see
that there are continous maps $\mathcal{M}(x,y)\to \Omega
(T(L,L'))$. We may then use these maps to construct a complex with
coefficients in $C_{\ast}(\Omega (T(L,L'))$. Clearly, there is an
obvious map $T(L,L')\to L$ and it is precisely this map which,
after looping, changes the coefficients  of this complex into
those of the extended Floer complex.

b. At this point it is worth mentioning why using representing
chain systems is useful in our constructions. Indeed, for the
extended Morse complex representing chain systems are not really
essential: the moduli spaces $M^{x}_{y}$ are triangulable in a way
compatible with the boundary formula and so, to represent this
moduli space inside the loop space $\Omega M$, we could use
instead of the chain $a^{x}_{y}$ a chain given by the sum of the
top dimensional simplexes in such a triangulation. This is
obviously, a simpler and more natural approach but it has the
disadvantage that it does not extend directly to the Floer case.
The reason is that it is not known whether the Floer moduli spaces
$\mathcal{M}(x,y)$ admit coherent triangulations (even if this is
likely to be the case). } \end{rem}

The chain complex $\mathcal{C}(L,L';H)$ admits a natural degree
filtration which is given by
\begin{equation}\label{eq:degree_filtr}F^{k}\mathcal{C}(L,L';H,J)=
C_{\ast}(L)\otimes \Z/2< x\in I(L,L';\eta,H) : |x|\leq k>~.~
\end{equation}

It is clear that this filtration is differential. Therefore, there
is an induced spectral sequence which will be denoted by
$\mathcal{E}(L,L'; H,J)=(\mathcal{E}^{r}_{p,q}, D^{r})$. We write
$\mathcal{E}(L,L';J)=\mathcal{E}(L,L';0,J)$. For convenience we
have omitted $\eta$ from this notation (the relevant components of
the paths spaces $\mathcal{P}(L,L')$ will be clear below).

Here is the main result concerning this spectral sequence.

\begin{theo}\label{theo:strips_ss}
For any two regular pairs $(H,J),(H',J')$, the spectral sequences
$\mathcal{E}(L,L';H, J)$ and $\mathcal{E}(L,L'; H',J')$ are
isomorphic up to translation for $r\geq 2$. Moreover, if $\phi$ is
a hamiltonian diffeomorphism , then $\mathcal{E}(L,L';J)$ and
$\mathcal{E}(L,\phi(L');J')$ are also isomorphic up to translation
for $r\geq 2$ (whenever defined). The second term of this spectral
sequence is $\mathcal{E}^{2}(L,L';H,J)\approx H_{\ast}(\Omega
L)\otimes HF_{\ast}(L,L')$ where $HF_{\ast}(-,-)$ is the Floer
homology. Finally, if $L$ and $L'$ are hamiltonian isotopic, then
$\mathcal{E}(L,L';J)$ is isomorphic up to translation to the Serre
spectral sequence of the path-loop fibration $\Omega L\to PL\to
L$.
\end{theo}

Isomorphism up to translation of two spectral sequences
$E^{r}_{p,q}$, $F^{r}_{p,q}$ means that there exists a $k\in \Z$
and chain isomorphisms $\phi^{r}: E^{r}_{p,q}\to F^{r}_{p+k,q}$.
This notion appears naturally here because the choice of the
element $x_{0}\in I(L,L';H)$ with $|x_{0}|=0$ is arbitrary.
 A different choice will simply lead to a translated spectral
sequence.
As follows from the discussion in \S\ref{susubsec:appli} B, it is possible
to replace this choice of grading with one that only depends on the path $\eta$.
However, this might make the absolute degrees fractionary and, as the choice of
$\eta$ is not canonical, the resulting spectral sequence will still be invariant only
up to translation.

\

The outline of the proof of this theorem is as follows (see
\cite{BaCo} for details). First, in view of the naturality
properties of Floer's construction, it is easy to see that the
second invariance claim in the statement is implied by the first
one. Now, we consider a homotopy $G$ between $H$ and $H'$ as well
as a one-parameter family of almost complex structures $\bar{J}$
relating $J$ to $J'$. For $x\in I(L,L';H,J)$ and $y\in
I(L,L';H',J')$ we define moduli spaces $\mathcal{N}(x,y)$ which
are solutions of an equation similar to (\ref{eq:floer}) but
replaces $H$ with $G$, $J$ with $\bar{J}$ (and takes into account
the additional parameter - this is a standard construction in
Floer theory). These moduli spaces have properties similar to the
$\mathcal{M}(x,y)$'s. In particular they admit compactifications
which are manifolds with boundary so that the following formula is
valid
$$\partial\overline{\mathcal{N}}(x,y)=\bigcup_{z\in I(L,L';H)}\overline{\mathcal{M}}
(x,z)\times \overline{\mathcal{N}}(z,y)\cup \bigcup_{z'\in
I(L,L';H')}\overline{\mathcal{N}}(x,z')\times\overline{\mathcal{M}}(z',y)~.~$$
The representing chain idea can again be used in this context and
it leads to coefficients $n^{x}_{y}\in C_{\ast}(\Omega L)$. If we
group these coefficients in a matrix $\mathcal{B}$ and we group
the coefficients of the differential of $\mathcal{C}(L,L';H,J)$ in
a matrix $\mathcal{A}$ and the coefficients of
$\mathcal{C}(L,L';H',J')$ in a matrix $\mathcal{A}'$, then the
relation above implies that we have:
\begin{equation}\label{eq:morph}
\partial \mathcal{B}=\mathcal{A}\cdot\mathcal{B}+\mathcal{B}\cdot\mathcal{A}'~.~
\end{equation}
It follows that the module morphism
$$\phi_{G,\bar{J}}:\mathcal{C}(L,L';H,J)\to \mathcal{C}(L,L';H,J)$$
which is the unique extension of
$$\phi_{G,\bar{J}}(x)=\sum_{y} r^{x}_{y}\otimes y,\ \forall
x\in I(L,L';H)$$ is a chain morphism.  Moreover, the chain
morphism constructed above preserves filtrations (of course, to
for this it is required that the choices for the point $x_{0}$
with $|x_{0}|=0$ for our two sets of data be coherent - this is
why the isomorphisms are ``up to translation"). After verifying
that $\mathcal{E}^{2}\approx H_{\ast}(\Omega L)\otimes
FH_{\ast}(L)$ for both spectral sequences it is not difficult to
see that $\phi_{G}$ induces an isomorphism at the
$\mathcal{E}^{2}$-level of these spectral sequences \cite{BaCo}.
Hence it also induces an isomorphism for $r >2$.

For the last point of the theorem we use Floer's reduction of the
moduli spaces $\mathcal{M}_{J,L,L'}(x,y)$ of pseudoholomorphic
strips to moduli spaces of Morse flow lines $M^{x}_{y}(f)$. In
short, this shows \cite{Fl1} \cite{Fl2} that for certain choices
of $J$, $f$ and $L'$ which is hamiltonian isotopic to $L$ we have
homeomorphisms $\psi_{x,y}:\mathcal{M}(x,y)\to M^{x}_{y}$ which
are compatible with the compactification and with the boundary
formulae. This means that with these choices we have an
isomorphism $\mathcal{C}(L,L')\to \mathcal{C}(f)$ and it is easy
to see that this preserves the filtrations of these two chain
complexes. By Theorem \ref{theo:serre_ss} this completes the
outline of proof.

\begin{rem}{\rm It is shown in \cite{BaCo} that
 the $\mathcal{E}^{1}$-term
of this spectral sequence has also some interesting invariance
properties. }
\end{rem}

\subsubsection{Applications.}\label{susubsec:appli}
We will discuss here a number of direct corollaries of Theorem
\ref{theo:strips_ss} most (but not all )of which appear in
\cite{BaCo}.

\

\underline{A. Localization and Hofer's metric.} An immediate
adaptation of Theorem \ref{theo:serre_ss} provides a statement
which is much more flexible. This is a ``change of coefficients"
or ``localization" phenomenon that we now describe.

Assume that $f:L\to X$ is a smooth map. Then we can consider the
induced map $\Omega f :\Omega L\to \Omega X$ and we may use this
map to change the coefficients of $\mathcal{C}(L,L';H,J)$ thus
getting a new complex $$\mathcal{C}_{X}(L,L';H,J)=(C_{\ast}(\Omega
X)\otimes \Z/2<I(L,L';H)>, D_{X})$$ so that
$D_{X}(x)=\sum_{y}(\Omega f)_{\ast}(m^{x}_{y})\otimes y$ where
$m^{x}_{y}$ are the coefficients in the differential $D$
 of $\mathcal{C}(L,L';H,J)$ (compare with (\ref{eq:differential_coeff})).
This complex behaves very much like the one studied in Theorem
\ref{theo:strips_ss}. In particular, this  complex admits a
similar filtration and the resulting spectral sequence,
$\mathcal{E}_{X}(L,L')$, has the same invariance properties as
those in the theorem and, moreover, for $L$, $L'$ hamiltonian
isoptopic this spectral sequence coincides with the Serre spectral
sequence of the fibration $\Omega X\to E\to L$ which is obtained
from the path-loop fibration $\Omega X\to PX\to X$ by pull-back
over the map $f$.

In particular, the homology of this complex coincides with the
singular homology of $E$. If $X$ is just one point, $\odot$, it is
easy to see that the complex $\mathcal{C}_{\odot}(L,L';H,J)$
coincides with the Floer complex.

The complex $\mathcal{C}_{X}(L,L';H,J)$ may also be viewed as a
sort of localization in the following sense. Assume that we are
interested to see what pseudo-holomorphic strips pass through a
region $A\subset L$. Then we may consider the closure $C$ of the
complement of this region and the space $L/C$ obtained by
contracting $C$ to a point. There is the obvious projection map
$L\to L/C$ which can be used in place of $f$ above. Now, if some
non-vanishing differentials appear in $\mathcal{E}_{L/C}(L,L')$
for $r\geq 2$, then it means that there are some coefficients
$m^{x}_{y}$ so that $|(m^{x}_{y})|>0$ and $(\Omega
f)_{\ast}(m^{x}_{y})\not=0$. This means that the map
$$\mathcal{M}(x,y)\stackrel{l^{x}_{y}} {\longrightarrow}\Omega
L\stackrel{\Omega f}{\longrightarrow}\Omega (C/L)$$ carries the
representing chain of $\mathcal{M}(x,y)$ to a nonvanishing chain
in $C_{\ast}(\Omega (L/C))$. But this means that the intersection
$\mathcal{M}'(x,y)\cap A$ is of dimension equal to $\mu(x,y)$.

The typical choice of region $A$ is a tubular neighbourhood of
some submanifold $V\hookrightarrow L$. In that case $L/C$ is the
associated Thom space.

 Let
$$\nabla(L,L')=\inf_{H, \phi^{H}_{1}(L)=L'}
 (\max_{x,t}H(x,t)-\min_{x,t}H(x,t))~.~$$
 be the Hofer distance between Lagrangians. It has been shown to
 be non-degenerate by Chekanov
 \cite{Chek} for symplectic manifolds with geometry bounded at
 infinity.

\begin{cor}\label{cor:strips_Maslov} Let $a\in H_{k}(L)$ be a non-trivial
homology class. If a closed submmanifold $V\hookrightarrow L$
represents the class $a$, then for any generic $J$ and any $L'$
hamiltonian isotopic to $L$, there exists a pseudoholomorphic
strip $u$ of Maslov index at most $n-k$ which passes through $V$
and which verifies:
 $$\int u^{\ast}\omega \leq \nabla (L,L')~.~$$
\end{cor}

In view of the discussion above, the proof is simple (we are using
here a variant of that used in \cite{BaCo}).

We start with a simple topological remark. Take $A$ to be a
tubular neighbourhood of $V$. Then $L/C=TV$ is the associated Thom
space. In the Serre spectral sequence of $\Omega (TV)\to P(TV)\to
TV$ we have that $d^{n-k}(\tau)\not=0$ where $\tau\in H_{n-k}(TV)$
is the Thom class of $V$. By Poincar\'e duality, there is a class
$b\in H_{n-k}(TV)$ which is taken to $\tau$ by the projection map
$L\to TV$ ($b$ corresponds to the Poincar\'e dual $a^{\ast}$ of
$a$ via the isomorphism $H^{n-k}(L)\approx H_{k}(L)$). This means
that $D^{n-k}(b)$ is not zero in $\mathcal{E}_{TV}(L,L';H)$ (for
any Hamiltonian $H$).

To proceed with the proof notice that, by the naturality
properties of the Floer moduli spaces, it is sufficient to show
that for any Hamiltonian $H$ (and any generic $J$) so that
$\phi^{H}_{1}(L)=L'$ there exists an element $u\in
\mathcal{M}'_{L,L,J,H}$ which is of Maslov index $(n-k)$ and so
that $u(\R,0)\cap
 A\not=\emptyset$ and $E_{L,L,H}(u)\leq ||H||$
 where $||H||=\max_{x,t}H(x,t)-\min_{x,t}H(x,t)$. We may assume
 that $\min_{x,t}H(x,t)=0$ and we let $K=\max_{x,t}H(x,t)$.
 We consider a Morse function $f:L\to \R$ which is very small in
 $C^{2}$ norm and we extend it to a function (also denoted by $f$)
 which is defined on $M$ and remains $C^{2}$ small. In
 particular, we
 suppose $\min_{x}f(x)=0$ and $\max_{x} f(x)<\epsilon$. We denote $\underline{f}=f-\epsilon$
 and $\overline{f}=f+K$. It follows that we may construct monotone
 homotopies $G:\overline{f}\simeq H$ and
 $G':H\simeq\underline{f}$. Consider the following action
 filtration of $\mathcal{C}_{TV}(L,L;H)$
 $$F_{v}\mathcal{C}_{TV}(L,L;J,H)=C_{\ast}(\Omega TV)\otimes \Z/2<x\in I(L,L;H):
 \mathcal{A}_{L,L,H}(x)\leq v>$$
 and similarly on the complexes $\mathcal{C}_{TV}(L,L;\overline{f})$
 and $\mathcal{C}_{TV}(L,L;\underline{f})$.
 It is obvious that this is a differential filtration and, if the
 choice of path $\eta$ (used to define the action functional, see
 (\ref{eq:action})) is the same for all the three hamiltonians
 involved, these monotone homotopies preserve these filtrations.

 We now denote
 $$\mathcal{C}=F_{K+\epsilon}\mathcal{C}_{TV}
 (L,L';\overline{f})/F_{-\epsilon}\mathcal{C}_{TV}(L,L';\overline{f}),$$

 $$\mathcal{C}'=F_{K+\epsilon}\mathcal{C}_{TV}(L,L';H)
 /F_{-\epsilon}\mathcal{C}_{TV}(L,L';H),$$

$$\mathcal{C}''=F_{K+\epsilon}\mathcal{C}_{TV}(L,L';\underline{f})
 /F_{-\epsilon}\mathcal{C}_{TV}(L,L';\underline{f})$$

These three complexes inherit degree filtrations and there are
associated spectral sequences $\mathcal{E}(\mathcal{C})$,
$\mathcal{E}(\mathcal{C}')$, $\mathcal{E}(\mathcal{C}'')$.
 We have induced morphisms $\phi_{G}:\mathcal{C}\to \mathcal{C}'$
 and $\phi_{G'}:\mathcal{C}'\to \mathcal{C}''$ which also induce
 morphisms among these spectral sequences. Moreover, as
 $\mathcal{C}$ and $\mathcal{C''}$ both coincide with
 $\mathcal{C}_{TV}(f)$ (because $f$ is very $C^{2}$-small and $0\leq f(x)<\epsilon$),
 the composition $\phi_{G'}\circ \phi_{G}$ induces an
 isomorphism of spectral sequences for $r\geq 2$
 (here $\mathcal{C}_{TV}(f)$ is the extended Morse complex obtained
 from $\mathcal{C}(f)$ by changing the coefficients by the map
 $L\to TV$).  But, as the class $b$ has the property that
 its $D^{n-k}$ differential is not trivial in
 $\mathcal{E}(\mathcal{C})$, this implies that
 $D^{n-k}(\phi_{G}(b))\not=0$ which is seen to immediately imply
 that there is some moduli space $\mathcal{M}'(x,y)$ of dimension
 $n-k$
with $\mathcal{A}_{L,L';H}(x),\mathcal{A}_{L,L';H}(y)\in
[-\epsilon,K+\epsilon]$ and $\mathcal{M'}(x,y)\cap
A\not=\emptyset$. Therefore there are $J$-strips passing through
$A$ which have Maslov index $n-k$ and area less than
$||H||+2\epsilon$. By letting $A$ tend to $V$ and $\epsilon\to 0$,
$||H||\to \nabla(L,L')$, these strips converge to strips with the
properties desired.

\

We may apply this even to $1\in H_{0}(L)$ and Corollary
\ref{cor:strips_Maslov} shows in this case that through each point
of $L$ passes a strip of Maslov index at most $n$ (again, for $J$
generic) and of area at most $\nabla(L,L')$. The case $V=pt$ was
discussed explicitly in \cite{BaCo}.

\begin{rem}\label{rem:area_esti}
{\rm  a. It is clear that the strips detected in this corollary
actually have a symplectic area which is no larger than
$c(b;H)-c(1;H)$ where $c(x;H)$ is the spectral value of the
homology class $x$ relative to $H$, $$c(x;H)=\inf\{v\in\R : x\in
Im(H_{\ast}(FC_{\leq v})(L,L;H)\to FH_{\ast}(L,L;H))\}$$
 where $FC_{\leq v}(L,L;H)$ is the Floer complex of $L,L,H$ generated
 by all the elements of $I(L,L;\eta, H)$ of action smaller or
 equal than $v$; $FH_{\ast}(L,L;H)$ is the Floer homology.
 Under our assumptions we have a canonical isomorphism
 (up to translation) between $HF_{\ast}(L,L, H)$ and $H_{\ast}(L)$
 so we may view $b\in HF_{\ast}(L,L;H)$.

 b. Clearly, in view of Gromov compactness our result also implies that
 for any $J$ (even non regular) and for any $L'$ hamiltonian isotopic to
 $L$ and for any $x\in L\backslash L'$ there exists a $J$-holomorphic
 strip passing through $x$ which has area less than $\nabla(L,L')$.
 This result (without the area estimate) also
 follows from independent work of Floer \cite{Fl} and Hofer
 \cite{Ho}. Another method has been mentioned to us by
 Dietmar Salamon. It is based
 on starting with disks with their boundary on $L$ and which
 are very close to be constant maps.  Therefore, an appropriate
 evaluation defined on the moduli space of these disks is of
 degree $1$. Each of these disks is made out of two
 semi-disks which are pseudo-holomorphic and which are joined by a
 short semi-tube verifying the non-homogenous Floer equation for
 some given Hamiltonian $H$. This middle region is then
 allowed to expand till, at some point, it will necessarily produce
 a semi-tube belonging to some $\mathcal{M}'_{H}(x,y)$.
It is also possible to use the pair of pants product to produce
Floer orbits joining the ``top and bottom classes'' \cite{Sch2}.
Still, having simultaneous area and Maslov index estimates appears
to be more difficult by methods different from ours. Of course,
detecting strips of lower Maslov index so that they meet some
fixed submanifold is harder yet.}\end{rem}

Corollary \ref{cor:strips_Maslov} has a nice geometric
application.

\begin{cor}\label{cor:balls}
 Assume that, as before, $L$ and $L'$ are hamiltonian isotopic.
For any symplectic embedding $e:(B_{r},\omega_{0})\to M$
 so that $e^{-1}(L)=\R^{n}\cap B(r)$ and $e(B_{r})\cap L'=\emptyset$
 we have $\pi r^{2}/2\leq \nabla(L,L')$.
\end{cor}

This is proven (see \cite{BaCo}) by using a  variant of the
standard isoperimetric inequality: a $J_{0}$-pseudoholomorphic
surface in the standard ball $(B_{r},\omega_{0})$ of radius $r$
whose boundary is on $\partial B_{r}\cup \R^{n}$ has area at least
$\pi r^{2}/2$.

Clearly, this implies the non-degeneracy result of Chekanov that
was mentioned before under the connectivity conditions that we
have always assumed till this point.

\

\underline{B. Relaxing the connectivity conditions.} We have
worked till now under the assumption that
\begin{equation}\label{eq:connect}
L, L' \ {\rm are\
simply-connected\ and\  } \omega|_{\pi_{2}(M)}=c_{1}|_{\pi_{2}(M)}=0
\end{equation}
These requirements were used in a few important places: in the
definition of the action functional, the definition of the Maslov
index, the boundary product formula (because they forbid
bubbling). Of these, only the bubbling isssue is in fact
essential: the boundary formula is precisely the reason why
$d^{2}=0$ as well as the cause of the invariance of the resulting
homology.

We proceed below to extend the corollaries and techniques
discussed above to the case when all the connectivity conditions
are dropped but we assume that $L$ and $L'$ are hamiltonian
isotopic and only work below the minimal energy that could produce
some bubbling (this is similar to the last section of \cite{BaCo}
but goes beyond the cases treated there).

First, for a time dependent almost complex structure $J_{t}$,
$t\in [0,1]$, we define $\delta_{L,L'}(J)$ as the infimmum of the
symplectic areas of the following three types of objects:
\begin{itemize}
\item[-] the $J_{t}$-pseudoholomorphic spheres in $M$ (for $t\in
[0,1]$). \item[-] the $J_{0}$-pseudoholomorphic disks with their
boundary on $L$. \item[-] the $J_{1}$-pseudoholomorphic disks with
their boundary on $L'$.
\end{itemize}

By Gromov compactness this number is strictly positive.

\

We will proceed with the construction in the case when $L=L'$ and
in the presence of a hamiltonian $H$. We shall assume that the
pair $(H,J)$ is regular in the sense that the moduli spaces of
strips defined below, $\mathcal{M}(x,y)$, are regular.

\

We take the fixed reference path $\eta$ to be a constant  point in
$L$ (see (\ref{eq:action})). Denote
$\mathcal{P}_{0}(L,L)=\mathcal{P}_{\eta}(L,L)$ and consider in
this space the base point given by $\eta$. Notice that there is a
morphism $\omega : \pi_{1}\mathcal{P}_{0}(L,L)\to \R$ obtained by
integrating $\omega$ over the disk represented by the element
$z\in\pi_{1}\mathcal{P}_{0}(L,L)$ (such an element can be viewed
as a disk with boundary in $L$). Similarly, let
$\mu:\pi_{1}\mathcal{P}_{0}(L,L)\to \Z$ be the Maslov morphism.
Let $\mathcal{K}$ be the kernel of the morphism
$$\omega\times\mu:\pi_{1}\mathcal{P}_{0}(L,L)\to \R\times \Z~.~$$

The group $$\pi=\pi_{1}(\mathcal{P}_{0}(L,L))/\mathcal{K}$$ is an
abelian group (as it is a subgroup of $\R\times \Z$) and is of
finite rank. Let's also notice that this group is the quotient of
$\pi_{2}(M,L)$ by the equivalence relation $a\sim b$ iff
$\omega(a)=\omega(b), \mu(a)=\mu(b)$ (with this definition this
group is also known as the Novikov group). This is a simple
homotopical result. First $\mathcal{P}(L,L)$ is the homotopy
pull-back of the map $L\to M$ over the map $L\to M$. But this
means that we have a fibre sequence $F\to \mathcal{P}(L,L)\to L$
with $F$ the homotopy fibre of $L\to M$ and that this fibre
sequence admits a canonical section. This implies that
$$\pi_{1}\mathcal{P}_{0}(L,L)\approx\pi_{1}(F)\times\pi_{1}(L)~.~$$
But $\pi_{1}(F)=\pi_{2}(M,L)$. As both $\omega$ and $\mu$ are
trivial on $\pi_{1}(L)$ the claim follows. It might not be clear
at first sight why $\mu$ is null on $\pi_{1}(L)$ here. The reason
is that the term $\pi_{1}(L)$ in the product above is the image of
the map induced in homotopy by $$j_{L}:L\to
\mathcal{P}_{0}(L,L)~.~$$ This map associates to a point in $L$
the constant path. Consider a loop $\gamma(s)$ in $L$. Then
$j_{L}\circ\gamma$ is a loop in $\mathcal{P}_{0}(L,L)$ which at
each moment $s$ is a constant path. We now need to view this loop
as the image of a disk and $\mu([j_{L}(\gamma)])$ is the Maslov
index of this disk. But this disk is null homotopic so
$\mu([j_{L}(\gamma)])=0$.

\

Consider the regular covering $p:\mathcal{P}'_{0}(L,L)\to
\mathcal{P}_{0}(L,L)$ which is associated to the group
$\mathcal{K}$. We fix an element
$\eta_{0}\in\mathcal{P}'_{0}(L,L)$ so that $p(\eta_{0})=\eta$.
Clearly, the action functional
$$\mathcal{A}'_{L,L,H}:\mathcal{P}'_{0}(L,L)\to \R$$
may be defined by essentially the same formula as in
(\ref{eq:action}):
$$ \mathcal{A}'_{L,L,H}(x)=
-\int (p\circ u)^{\ast}\omega +\int_{0}^{1}H(t,(p\circ x)(t))dt
$$
where $u:[0,1]\to \mathcal{P}'_{0}(L,L)$ is such that
$u(0)=\eta_{0}$, $u(1)=x$ and is now well-defined.

Let $I'(L,L,H)=p^{-1}(I(L,L,H))$. For $x,y\in I'(L,L,H)$ we may
define $\mu(x,y)=\mu(p\circ u)$ where $u:[0,1]\to
\mathcal{P}'_{0}(L,L)$ is a path that joins $x$ to $y$. This is
again well defined. For each $x\in I'(L,L,H)$ we consider a path
$v_{x}:[0,1]\to \mathcal{P}'_{0}(L,L)$ so that $v_{x}(0)=\eta_{0}$
and $v_{x}(1)=x$. The composition $p\circ v_{x}$ can be viewed as
a ``semi-disk" whose boundary is resting on the orbit $p(x)$ and
on $L$. Therefore, we may associate to it a Maslov index
$\mu(v_{x})$ \cite{RoSa} and it is easy to see that this only
depends on $x$. Thus we define $\mu(x)=\mu(v_{x})$ and we have
$\mu(x,y)=\mu(x)-\mu(y)$ for all $x,y\in I'(L,L,H)$.

\

To summarize what has been done till now: once the choices of
$\eta$ and $\eta_{0}$ are made, both the action functional
$\mathcal{A}':\mathcal{P}'_{0}(L,L,H)\to \R$ and the ``absolute"
Maslov index $\mu(-):I'(L,L,H)\to \Z$ are well-defined.

\

Fix an almost complex structure $J$. Consider two elements $x,y\in
I'(L,L,H)$. We may consider the moduli space which consists of all
paths $u:\R\to \mathcal{P}'_{0}(L,L)$ which join $x$ to $y$ and
are so that $p\circ u$ satisfies Floer's equation (\ref{eq:floer})
modulo the $\R$-action. If regularity is achieved, the dimension
of this moduli space is precisely $\mu(x,y)-1$. The action
functional $\mathcal{A}'$ decreases along such a solution $u$ and
the energy of $u$ (which is defined as the energy of $p\circ u$)
verifies, as in the standard case,
$E(u)=\mathcal{A}'(x)-\mathcal{A}'(y)$. Bubbling might of course
be present in the compactification of these moduli spaces. As we
only intend to work below the minimal bubbling energy
$\delta_{L,L}(J)$ we {\em artificially} put:
$$\mathcal{M}(x,y)=\emptyset\ {\rm if}\
\mathcal{A}'(x)-\mathcal{A}'(y)\geq \delta_{L,L}(J)$$ and, of
course, for $\mathcal{A}'(x)-\mathcal{A}'(y)<\delta_{L,L}(J)$,
$\mathcal{M}(x,y)$ consists of the elements $u$ mentioned above.
We only require these moduli spaces to be regular.

With this convention, for all $x,y\in I'(L,L,H)$ so that
$\mathcal{M}(x,y)$ is not void we have the usual boundary formula
(\ref{eq:Grom_comp}). Notice at the same time that this formula is
false for general pairs $x,y$ (and so there is no way to define a
Floer type complex at this stage).

\

Now consider a map $f:L\to X$ so that $X$ is simply-connected (the
only reason to require this is to insure that the Serre spectral
sequence does not require local coefficients).

We consider the group:
$$\mathcal{C}(L,L,H;X)=C_{\ast}(\Omega X)\otimes \Z/2<I'(L,L,H)>~.~$$
For $w\geq v\in \R$, we denote  $I'_{v,w}=\{x\in I'(L,L,H)\ :\
w\geq\mathcal{A}'(x)\geq v\}$ and we define the subgroup
$$\mathcal{C}_{w,v}(L,L,H;X)=C_{\ast}(\Omega X)\otimes \Z/2<I'_{v,w}>~.~$$

Suppose that $w-v \leq \delta_{L,L}(J)-\epsilon$. We claim that in
this case we may define a differential on
$\mathcal{C}_{v,w}(L,L,H;X)$ by the usual procedure. Consider
representing chain systems for all the moduli spaces
$\mathcal{M}(x,y)$ and let the image of these chains inside
$C_{\ast}(\Omega X)$ be respectively $\bar{m}^{x}_{y}$. Let $D$ be
the linear extension of the map given by
$$Dx=\sum_{y\in I'_{v,w}}\bar{m}^{x}_{y}\otimes y~.~$$

\begin{prop}\label{prop:extension}
The linear map $D$ is a differential. A generic monotone homotopy
$G$ between two hamiltonians $H$ and $H'$
$$\phi_{G}^{X}:\mathcal{C}_{v,w}(L,L,H, J;X)
\to\mathcal{C}_{v,w}(L,L,H',J;X)~.~$$ A monotone homotopy between
monotone homotopies $G$ and $G'$ induces a chain homotopy between
$\phi_{G}^{X}$ and $\phi_{G'}^{X}$ so that
$H_{\ast}(\phi_{G}^{X})= H_{\ast}(\phi_{G'}^{X})$.
\end{prop}

Now, $D^{2}x=\sum_{z}(\sum_{y}\bar{m}^{x}_{y}\cdot
\bar{m}^{y}_{z}+\partial \bar{m}^{x}_{z})\otimes z~.~$ In this
formula we have $\mathcal{A}'(x)-\mathcal{A}'(z)\leq
\delta_{L,L}(J)$ and, because the usual boundary formula
(\ref{eq:Grom_comp}) is valid in this range, all the terms vanish
so that $D^{2}(x)=0$. The same idea may be applied to a monotone
homotopy as well as to a monotone homotopy between monotone
homotopies and it implies the claim.

\begin{rem}\label{rem:point} {\rm a. If we take for the space $X$ a
single point $\ast$ we get a chain complex whose differential only
takes into account the $0$-dimensional moduli spaces and which is a
truncated version of Floer homology.

b. The complex $\mathcal{C}_{v,w}(L,L,H,J;X)$ admits a degree
filtration which is perfectly similar to the one given by
(\ref{eq:degree_filtr}). Let $\mathcal{E}\mathcal{C}_{v,w}(L,L,H,
J;X)$ be the resulting spectral sequence. Then, under the
restrictions in the Proposition \ref{prop:extension}, a monotone
homotopy $G$ induces a morphism of spectral sequences
$\mathcal{E}_{X}(\phi_{G})$ and two such homotopies $G$, $G'$
which are monotonously homotopic have the property that they
induce the same morphism
$\mathcal{E}^{r}_{X}(\phi_{G})=\mathcal{E}^{r}_{X}(\phi_{G'})$ for
$r\geq 2$. This last fact follows from Proposition
\ref{prop:extension} by computing
$\mathcal{E}^{2}_{X}(\phi_{G})=id_{H_{\ast}(\Omega X)}\otimes
H_{\ast}(\phi_{G}^{\ast})=\mathcal{E}^{2}_{X}(\phi_{G'})$ where
$\phi_{G}^{\ast}:\mathcal{C}_{v,w}(L,L,H;\ast)\to
\mathcal{C}_{v,w}(L,L,H';\ast)$.}
\end{rem}

\

Naturally, the next step is to compare our construction with its
Morse theoretical analogue.
Consider the map $j_{L}:L\to \mathcal{P}_{0}(L,L)$ and consider
$p:\tilde{L}\to L$ the regular covering obtained by pull-back from
$\mathcal{P}'_{0}(L,L)\to \mathcal{P}_{0}(L,L)$. Notice that,
because both compositions $\omega\circ\pi_{1}(j_{L})$ and
$\mu\circ\pi_{1}(j_{L})$ are trivial, it follows that the covering
$\tilde{L}\to L$ is trivial. Let $\bar{f}:\tilde{L}\to \R$ be
defined by $\bar{f}=f\circ p$ and consider $\mathcal{C}(\bar{f};X)$ the extended
Morse complex of $\bar{f}$ with coefficients changed by the map $\Omega
\tilde{L}\to \Omega X$. Notice that, in
general, the group $\pi$ acts on $I'(L,L,H)$ and we have the
formula:
$$\mathcal{A}'(gy)=\omega(g)+\mathcal{A}'(y),\ \mu(gy)=\mu(g)+\mu(y),
\forall y\in I'(L,L,H), \forall g\in\Pi~.~$$ In our particular
case, when $H=f$, we have $I'(L,L,H)=Crit(\bar{f})$. For each
point $x\in Crit(f)$ let $\bar{x}\in Crit(\bar{f})$ be the element
of $p^{-1}(x)$ which belongs to the component of $\tilde{L}$ which
also contains $\eta_{0}$. We then have
$\mathcal{A}'(\bar{x})=f(x)$ and $\mu(\bar{x})=ind_{\bar{f}}(\bar{x})$. The extended Morse complex
$\mathcal{C}(\bar{f};X)$ is therefore isomorphic to
$\mathcal{C}(f;X)\otimes \Z[\pi]$ and the action filtration is
determined by writting $\mathcal{A}'(x\otimes g)=f(x)+\omega(g)$.
The degree filtration induces, as usual, a spectral sequence
which will be denoted by $\mathcal{E}\mathcal{C}(f;X)$.

The remarks above together with
Theorem \ref{theo:strips_ss}  show that this spectral sequence consists of copies
the Serre spectral sequence of $\Omega X\to
E\to L$: one copy for each connected component of
$\tilde{L}$. We denote by $\mathcal{C}_{0}(f;X)$ and $\mathcal{E}\mathcal{C}_{0}(f;X)$
the copies of the extended complex and of the spectral sequence that correspond to the
connected component $L_{0}$ of $\tilde{L}$ which contains $\eta_{0}$.

\begin{prop}\label{prop:comp_morse} Suppose $||H||_{0}< \delta_{J}(L,L)$.
There exists a chain morphisms
$\phi:\mathcal{C}_{0}(f;X)\to \mathcal{C}_{0, ||H||}(L,L,H,J;X)$
and $\psi :\mathcal{C}_{0,||H||}(L,L,H,J;X) \to \mathcal{C}_{0}(f;X)$
which preserve the respective degree filtrations and so that $\psi\circ\phi$ induces
 an isomorphism at the $E^{2}$ level of the respective spectral sequences.
\end{prop}

To prove this proposition we shall use a different method
than the one used in Corollary \ref{cor:strips_Maslov}.  The comparison maps
$\phi$, $\psi$ will be constructed by the method introduced in \cite{PiSaSc}
and later used in \cite{Sch1} and \cite{Sch2}. Compared to
Proposition \ref{prop:extension} this is particularly efficient
because it avoids the need to control the bubbling thershold
along deformations of $J$.

\

We fix as before the Morse function $f$ as well as the pair $H,J$. To simplify
the notation we shall assume that $\inf H(x,t)=0$.
The construction of $\phi$ is based on defining certain moduli spaces
$\mathcal{W}(x,y)$ with $x\in Crit(\bar{f})$ and $y\in I'(L,L,H)$.
They consist of pairs $(u,\gamma)$ where $u:\R\to \mathcal{P}'_{0}(L,L)$,
$\gamma:(-\infty,0]\to \tilde{L}$ and if we put $u'=p(u)$, $\gamma'=p(\gamma)$
($p:\mathcal{P}'_{0}(L,L)\to \mathcal{P}_{0}(L,L)$ is the covering projection)
then we have:
$$u'(R\times\{0,1\})\subset L\ , \ \partial_{s}(u')+ J(u')\partial_{t}(u')+\beta(s)\nabla
H(u',t)=0\ , u(+\infty)=y$$ and
$$
\frac{d\gamma'}{dt}=-\nabla_{g}f(\gamma')\ ,\  \gamma(-\infty)=x \ , \
\gamma(0)=u(-\infty)~.~$$
 Here $g$ is a Riemannian metric so that $(f,g)$ is Morse-Smale and $\beta$ is a smooth
 cutt-off function which is increasing and vanishes for $s\leq 1/2$ and equals
 $1$ for $s\geq 1$. It is useful to view an element $(u,\gamma)$ as before as a
 semi-tube connecting $x$ to $y$. Under usual regularity assumptions
 these moduli spaces are manifolds of dimension $\mu(x)-\mu(y)$.

The energy of such an element $(u,\gamma)$ is defined in the obvious way by
$E(u,\gamma)=\int ||\partial_{s}u'||^{2}ds dt$. A simple computation
shows that:
$$E(u,\gamma)=I(u) +
\int_{\R\times [0,1]} (u')^{\ast}\omega-\int_{0}^{1}H(y(t))dt$$
where $I(u)=\int_{\R\times [0,1]} \beta'(s)H(u'(s),t)dsdt$.
If $x\in L_{0}$, then the energy verifies
$$E(u,\gamma)=I(u) -\mathcal{A}'(y)\leq \sup (H) -\mathcal{A}'(y)~.~$$

As before, we only want to work here under the bubbling threshold and we are only
interested in the critical points $x\in L_{0}$ so we put
$\mathcal{W}(x,y)=\emptyset$ for all those pairs $(x,y)$ with either $y\in I'(L,L,H)$ so that
$\mathcal{A}'(y)\not\in [0,||H||]$ or with $x\not\in L_{0}$.
This means that there is no bubbling in our moduli spaces. Thus we may apply the usual
procedure: compactification, representing chain systems, representation in the loop space
(for this step we need to choose a convenient way to parametrize the paths represented
by the elements $(u,\gamma)$).  Notice that the boundary of
$\overline{\mathcal{W}}(x,y)$ is the union of two types of pieces
$\overline{M}^{x}_{z}\times \overline{\mathcal{W}}(z,y)$ and $\overline{\mathcal{W}}(x,z)\times
\overline{\mathcal{M}}(z,y)$.
We then define $\phi(x)=\sum w^{x}_{y}\otimes y$
where $w^{x}_{y}$ is a cubical chain representing the moduli space $\mathcal{W}(x,y)$.
This is a chain map as desired.

We now proceed to construct the map $\psi$. The construction is prefectly similar:
we define moduli spaces $\mathcal{W}'(y,x)$, $y\in I'(L,L,J)$, $x\in Crit(\bar{f})$  except
that the pairs $(u,\gamma)$ considered here, start as semi-tubes and end as flow lines of $f$.
The equation verified by $u$ is similar to the one before but istead of the cut-off function
$\beta$ we use the cut-off function $-\beta$. For $x\in L_{0}$ the energy estimate in this
case gives $E(u,\gamma)\leq \mathcal{A}'(y)$.
By the same method as above we define $\mathcal{W}'(y,x)$ to be void whenever $x\not\in L_{0}$ or
$\mathcal{A}'(y)>||H||$ and we define
$\psi(y)=\sum \bar{w}^{y}_{x}\otimes x$ where $\bar{w}^{y}_{x}$ is a cubical chain
representing the moduli space $W'(y,x)$. Notice that because $E(u,\gamma)\leq \mathcal{A}'(y)$
this map $\psi$ does in fact vanish on $\mathcal{C}_{0}(L,L,H,J;X)$ and so it induces
a chain map (also denoted by $\psi$) as desired.

The next step is to notice that the composition $\psi\circ\phi$ induces an isomorphism
at $E^{2}$. This is equivalent to showing that $H_{\ast}(\psi\circ\phi)$ is an isomorphism
for $X=\ast$. In turn, this fact follows by now standard deformation
arguments as in \cite{PiSaSc}.

\begin{cor}\label{cor:general_strips}
Assume that $L$ and $L'$ are hamiltonian isotopic and suppose that
$J$ is generic. If $\nabla(L,L')< \delta_{L,L'}(J)$, then the
statement of Corollary \ref{cor:strips_Maslov} remains true (for
$J$) without the connectivity assumption (\ref{eq:connect}).
\end{cor}

Notice that if $H$ is a hamiltonian so that $\phi^{H}_{1}(L)=L'$
and $J_{H}=(\phi^{H})_{\ast}(J)$ then, by the naturality described
in \S\ref{subsec:strips}, we have:
$$\delta_{L,L}(J_{H})=\delta_{L,L'}(J)~.~$$
This implies, again by this same naturality argument, that the
problem reduces to finding appropriate semi-tubes whose detection
comes down to showing the non-vanishing of certain differentials
in $\mathcal{E}\mathcal{C}_{0,w}(L,L,H;TV)$ for some well chosen
 $w<\delta_{L,L}(J_{H})$. But this immediatly follows from Proposition
 \ref{prop:comp_morse} by the same topological argument as the one used
 in the proof of Corollary \ref{cor:strips_Maslov}

\

We formulate the geometric consequence which corresponds to
Corollary \ref{cor:balls}. For two lagrangians $L$ and $L'$ the
following number has been introduced in \cite{BaCo}: $B(L,L')$ is
the supremum of the numbers $r\geq 0$ so that there exists a
symplectic embedding
$$ e: (B(r),\omega_{0})\to (M,\omega)$$ so that
$e^{-1}(L)= \R^{n}\cap B(r)$ and $Im(e)\cap L'=\emptyset$.

\begin{cor}\label{cor:ball_bubble} There exists an almost
complex structure $J$ so that we have the inequality:
$$\nabla(L,L')\geq \min\{\delta_{L,L'}(J)\ ,\ \frac{\pi}{2}
B(L,L')^{2}\}~.~$$
\end{cor}

 Clearly, this implies that $\nabla(-,-)$ is non-degenerate
 in full generality (and recovers, in particular, the fact that
 the usual Hofer norm for Hamiltonians is non-dgenerate). It
 is useful to also notice as in \cite{BaCo} that this result
 is a Lagrangian version of the
 usual capacity - displacement energy inequality \cite{LaMc}.
 Indeed, this inequality (with the factor $\frac{1}{2}$)
 is implied by the following statement which has been
 conjectured  to hold for any two compact lagrangians in a symplectic
 manifold \cite{BaCo}:

\begin{equation}\label{eq:conj} \nabla(L,L')\geq
\frac{\pi}{2}B(L,L')^{2}~.~
\end{equation}

This remains open. An even stronger conjecture is the following:

\

\begin{conj}\label{conj:main} For any two hamiltonian isotopic closed lagrangians $L,L'\subset
(M,\omega)$ and for any almost complex structure $J$ which compatible with
$\omega$ and any point $x\in L\backslash L'$ there exists a
pseudoholomorphic curve $u$ which is either a strip resting on $L$
and on $L'$ or a pseudoholomorphic disk with boundary in $L$ so
that $x\in Im(u)$ and $\int u^{\ast}\omega\leq \nabla(L,L')$.
\end{conj}
\

By the isoperimetric inequality used earlier in this paper, it
follows that this statement implies (\ref{eq:conj}). There is a
substantial amount of evidence in favor of this conjecture:
\begin{itemize}
\item[-] as explained in this paper, in the absence of any
pseudoholomorphic disks (that is, when $\omega|_{\pi_{2}(M,L)}=0$)
it was proven in \cite{BaCo}. \item[-]the statement in Corollary
\ref{cor:ball_bubble} shows that the area estimate is not
unreasonable. \item[-]one striking consequence of Conjecture
\ref{conj:main} is that if the disjunction energy of the
lagrangian $L$ is equal to $E_{0}<\infty$, then, for any $J$ as in
the statement and any $x\in L$ there is a pseudoholomorpic disk of
area at most $E_{0}$ which passes through $x$. When $L$ is
relatively spin, this is indeed true and follows from recent work
of the second author joint with Fran\c{c}ois Lalonde. By the same
geometric argument as above we deduce a nice consequence. Define
the relative (or \emph{real}) Gromov radius of $L$, $Gr(L)$, to be
the supremum of the positive numbers $r$ so that there exists  a
symplectic embedding $e:(B(r),\omega_{0})\to (M,\omega)$ with the
property that $e^{-1}(L)= \R^{n}\cap B(r)$, then $\pi
(Gr(L))^{2}/2\leq E_{0}$ (where $E_{0}$, as before, is the
disjunction energy of $L$). It is also useful to note that if $L$
is the zero section of a cotangent bundle, then $Gr(L)=\infty$.
\end{itemize}

There are numerous other interesting consequences of Conjecture
\ref{conj:main} besides (\ref{eq:conj}).  To conclude,
\ref{conj:main} appears to be a statement worth investigating.


\begin{thebibliography}{10}

\bibitem{BaCo} J.-F. Barraud, O. Cornea, \emph{Lagrangian
Intersections and the Serre spectral sequence}, Prepint, November
2003.

\bibitem{BeCo} M.Betz, R.Cohen \emph{Graph moduli spaces and
cohomology operations}, Turkish J. Math. 18 (1994), no. 1, 23--41.

\bibitem{Chek}
Y. Chekanov, \emph{Invariant Finsler metrics on the space of
Lagrangian embeddings}, Math. Z. 234,(2000),  605\u2013-619.

\bibitem{CoJoSe1}{R.L.Cohen, J.D.S.Jones, G.B.Segal}, {Morse theory and classifying
spaces}{}{}{Preprint}

\bibitem{CoJoSe2}{R.Cohen, J.D.S.Jones, G.Segal},{ Floer's infinite dimensional Morse
theory and homotopy theory}{The Floer memorial Volume,
Birkhauser}{}{ (1995)}

\bibitem{Co1}
O. Cornea, \emph{Homotopical Dynamics II: Hopf invariants,
smoothings and the Morse complex}, Ann. Scient. Ec. Norm. Sup.  35
(2002) 549--573.

\bibitem{Co2}
O. Cornea, \emph{Homotopical Dynamics IV: Hopf invariants and
Hamiltonian flows}, Communications on Pure and Applied Math.  55
(2002), 1033--1088.

\bibitem{Co3}
O.Cornea, \emph{New obstructions to the thickening of
CW-complexes}, Proc. Amer. Math. Soc., 132, (2004) no. 9,
2769--2781.

\bibitem{CLOT}
O. Cornea, G. Lupton, J. Oprea, D. Tanr\'e,
\emph{Lusternik-Schnirelmann category}, AMS Monographs in
Mathematics, 2003.

\bibitem{Fl}
A. Floer, \emph{Cuplength estimates on Lagrangian intersections},
Comm. Pure Appl. Math. 42 (1989), no. 4, 335--356.

\bibitem{Fl1}
A. Floer, \emph{ Morse theory for Lagrangian intersections}, J.
Differential Geom. 28 (1988), no. 3, 513--547.

\bibitem{Fl2}
A. Floer,  \emph{Witten's complex and infinite-dimensional Morse
theory}, J. Differential Geom. 30 (1989), no. 1, 207--221.


\bibitem{Fr}
J.Franks \emph{Morse-Smale flows and homotopy theory}, Topology
18, (1979), 199--215.

\bibitem{Fu}
K. Fukaya,  \emph{ $A\sp \infty$-category, and Floer homologies},
Proceedings of GARC Workshop on Geometry and Topology '93 (Seoul,
1993), 1--102, Lecture Notes Ser., 18, Seoul Nat. Univ., Seoul,
1993.

\bibitem{Ho}
H. Hofer, \emph{Lusternik-Schnirelman-theory for Lagrangian
intersections},
 Ann. Inst. H. Poincar\'e  Anal. Non Lin\'eaire 5 (1988), no 4, 465--499.

\bibitem{LaMc}
F. Lalonde, D. McDuff, \emph{The Geometry of Symplectic Energy}
Ann. of Math. (2) 141 (1995), no. 2, 349--371.


\bibitem{Mil}
J. Milnor, \emph{Lectures on the $h$-cobordism theorem}, Notes by
L. Siebenmann and J. Sondow Princeton University Press, Princeton,
N.J. 1965.

\bibitem{Mil2}{J. Milnor},\emph{Topology from the differentiable
viewpoint}, Based on notes by David W. Weaver The University Press
of Virginia, Charlottesville, Va. 1965.

\bibitem{PiSaSc} S.Piunikhin, D. Salamon, M. Schwarz,
\emph{Symplectic Floer-Donaldson theory and quantum cohomology},
Contact and symplectic geometry (Cambridge, 1994), 171--200,
Publ. Newton Inst., 8, Cambridge Univ. Press, Cambridge, 1996.

\bibitem{PuRo}{Ch.C.Pugh, C.Robinson},
\emph{The $C^{1}$ closing Lemma, including
Hamiltonians},{Erg.Theory \& Dyn. Sys.}{3}{(1983) 261--313.}

\bibitem{RoSa}
J. Robbin, D. Salamon, \emph{Asymptotic behaviour of holomorphic
strips}. Ann. Inst. H. Poincaré Anal. Non Linéaire 18 (2001), no.
5, 573--612.

\bibitem{Sal1}{D.Salamon}, \emph{Connected simple systems and the Conley
index of isolated invariant sets}, {Trans. of the A.M.S.},{291},
{(1985) 1--41.}

\bibitem{Sal2}
D. Salamon, \emph{Lectures on Floer Homology}, Symplectic Geometry
and Topology, edited by Y. Eliashberg and L. Traynor, IAS/Park
City
 Mathematics series,  (1999), 143--230.

\bibitem{Sch}
M.Schwarz, \emph{Morse homology}, Progress in Mathematics, 111.
Birkh\"auser Verlag, Basel, 1993.

\bibitem{Sch1}
M. Schwarz, \emph{ A quantum cup-length estimate for symplectic
fixed points}, Invent. Math. 133 (1998), no. 2, 353--397.

\bibitem{Sch2}
M. Schwarz, \emph{On the action spectrum for closed symplectically
aspherical manifolds}, Pacific J. Math. 193 (2000), no. 2,
419--461.


\bibitem{Sm}
S.Smale, \emph{Differentiable dynamical systems}, Bull. A.M.S.,
73,  (1967), 747--817.

\bibitem{Viterbo}
C.Viterbo, \emph{Intersection de sous-vari\'et\'es lagrangiennes,
fonctionnelles d'action et indice des systèmes hamiltoniens},
Bull. Soc. Math. France 115 (1987), no. 3, 361--390.

\bibitem{We}{J.Weber}, The Morse-Witten complex via dynamical
systems, preprint Winter 2004.

\bibitem{Wit}{E.Witten}, {Supersymmetry and Morse theory}, {J. of Diff.
Geometry}, {17}, { (1982), 661--692.}

\end{thebibliography}
\end{document}